\documentclass[11pt]{article}

\usepackage{slashed}
\usepackage{amsmath}
\usepackage{amssymb}
\usepackage{amsthm}
\usepackage{latexsym}
\usepackage{color}
\usepackage{graphicx}
\usepackage{color}
\usepackage{hyperref}
\hypersetup{
    colorlinks,
    citecolor=blue,
    filecolor=blue,
    linkcolor=blue,
    urlcolor=blue
}

\DeclareSymbolFont{calletters}{OMS}{cmsy}{m}{n}
\DeclareSymbolFontAlphabet{\mathcal}{calletters}

\def\be{\begin{eqnarray}}
\def\ee{\end{eqnarray}}

\def\b*{\begin{eqnarray*}}
\def\e*{\end{eqnarray*}}

\newtheorem{Theorem}{Theorem}[section]

\newtheorem{Proposition}[Theorem]{Proposition}

\newtheorem{Assumption}[Theorem]{Assumption}

\newtheorem{Lemma}[Theorem]{Lemma}

\newtheorem{Remark}[Theorem]{Remark}
\newtheorem{Example}[Theorem]{Example}

\makeatletter \@addtoreset{equation}{section}

\usepackage{color}

\def \C{\mathbb{C}}

\def \E{\mathbb{E}}

\def \L{\mathbb{L}}
\def \P{\mathbb{P}}

\def \R{\mathbb{R}}

\def \N{\mathbb{N}}

\def\Fc{{\cal F}}

\def\Kc{{\cal K}}

\def\Sc{{\cal S}}

\def\Wc{{\cal W}}

\def \0{\mathbf{0}}

\def\1{{\bf 1}}
\def \proof{{\noindent \bf Proof.\quad}}

\def \no{\noindent}
\def \ep{\hbox{ }\hfill{ ${\cal t}$~\hspace{-5.1mm}~${\cal u}$}}

 \newcommand{\eps}{\varepsilon}

\newcommand{\bea}{\begin{eqnarray}}
\newcommand{\bes}{\begin{subequations}}
\newcommand{\ees}{\end{subequations}}
\newcommand{\bgt}{\begin{gather}}
\newcommand{\egt}{\begin{gather}}
\newcommand{\eea}{\end{eqnarray}}
\newcommand{\beaa}{\begin{eqnarray*}}
\newcommand{\eeaa}{\end{eqnarray*}}

\newcommand{\EE}{{\mathbb E}}
\newcommand{\RR}{{\mathbb R}}

\title{Branching diffusion representation for nonlinear Cauchy problems
and Monte Carlo approximation}

\author{Pierre Henry-Labord\`ere\thanks{Soci\'et\'e G\'en\'erale, Global Market Quantitative Research,
pierre.henry-labordere@sgcib.com}
        \and Nizar Touzi\thanks{Ecole Polytechnique Paris, Centre de Math\'ematiques Appliqu\'ees,
        nizar.touzi@polytechnique.edu. This work benefits from the financial support of the ERC Advanced Grant 321111, and the Chairs {\it Financial Risk} and {\it Finance and Sustainable Development}. }}

\date{\today}

\begin{document}

\maketitle

\begin{abstract}\no We provide a probabilistic representations of the solution of some semilinear hyperbolic and high-order PDEs based on branching diffusions. These representations pave the way for a Monte-Carlo approximation of the solution, thus bypassing the curse of dimensionality. We illustrate the numerical implications in the context of some popular PDEs in physics such as nonlinear Klein-Gordon equation, a simplified scalar version of the Yang-Mills equation, a fourth-order nonlinear beam equation and the Gross-Pitaevskii PDE as an example of nonlinear Schr\"{o}dinger equations.
\end{abstract}

\section{Introduction}

Similar to the intimate connection between the heat equation and the Brownian motion, linear (second-order) parabolic partial differential equations are connected to stochastic processes. More precisely, the Feynman-Kac formula states that a linear parabolic PDE with infinitesimal generator ${\cal L}:=b\cdot\partial_x+\frac12{\rm Tr}[\sigma\sigma^{\rm T}\partial^2_{xx^{\rm\footnotesize T}}]$, and terminal condition $f_1(\cdot)$ at time $T$, can be written as a conditional expectation of $f_1(X_T)$ involving the It\^o process $X$ associated with the generator $\cal L$. This connection allows to device numerically approximations of the solution of such PDEs by probabilistic (Monte Carlo) methods, which represents a clear advantage in high dimensional problems as the error estimate induced by the central limit theorem is dimension-free.

\no An important focus was put on the extension to nonlinear (second-order) parabolic PDEs, see e.g. \cite{guy} for a quick review of existing methods.

\no - A first attempt was achieved by exploiting the stochastic representation by means of backward stochastic differential equations, see Bally \& Pag\`es \cite{BallyPages}, Bouchard \& Touzi \cite{BouchardTouzi}, Zhang \cite{Zhang}, and the extension to the fully nonlinear setting by Fahim, Touzi \& Warin \cite{FTW}. The induced numerical method involves repeated computations of conditional expectations, resulting in a serious dimension dependence of the corresponding numerical methods. In fact, as highlighted in \cite{FTW}, this method can be viewed as part of the traditional finite-elements algorithm.

\no - More recently, Henry-Labord\`ere \cite{phl1} suggested a numerical method based on an extension of the McKean \cite{McKean} branching diffusion representation of the so-called KPP equation (Kolmogorov-Petrovskii-Piskunov) to a class of semilinear second order parabolic PDEs with power nonlinearity in the value function $u$. The resulting algorithm is purely probabilistic. In particular its convergence is controlled by the central limit theorem, and is independent of the dimension of the underlying state. The validity of this method in the path-dependent case, and for further analytic nonlinearities which are of the power type in the triple $(u,\partial_xu,\partial^2_{xx^{\rm T}}u)$ is analyzed in \cite{phl1,phl2,phl3,phl4}. We also refer to Agarwal and Claisse \cite{aga} for the extension to elliptic semilinear PDEs, and Bouchard, Tan, Warin \& Zou \cite{BouchardTanWarinZou} for Lipschitz nonlinearity in the pair $(v,\partial_xv)$. A critical ingredient for the extension is the use of Galton-Watson trees weighted by some Malliavin automatic differentiation weights.

\no Our objective in this  paper is to show that the above branching diffusion approach extends to more general Cauchy problems, including hyperbolic and higher order PDEs. Such an extension may seem to be unexpected due to the intimate connection between parabolic second order PDEs and diffusions generators. However, probabilistic representations of some specific examples of hyperbolic PDEs did appear in the previous literature. The first such relevant work traces back to Kac \cite{kac} in the context of the one-dimensional telegrapher equation:
 \b*
 \partial_{tt} u- c\,\partial_{xx}u+({2 \beta})\partial_t u
 =
 0,
 &u(0,x) \equiv f_1(x),&
 \partial_t u(0,x)\equiv f_2(x),
 \e*
where $\beta$ and $c$ are constant parameters, and the boundary data $f_1,f_2$ are some given functions. Observing that we may express $u(t,x)=\frac12[u_0+u_1](t,x)$, where
the pair $(u_0,u_1)$ solves the coupled system of first order PDEs:
 $$
 \partial_t u_i +(-1)^{i}c\, \partial_x u_i -\beta(u_{1-i}-u_i)
 =
 0,
 ~
 u_i(0,x)
 =
 f_1(x)-{(-1)^i \over c}\int_0^x f_2(y) dy,
 ~i=0,1,
 $$
\no A branching mechanism representation of $u_0$ and $u_1$ is obtained by following McKean's representation for (interacting) KPP equations (with zero diffusion).

\no We next mention the work of Dalang, Mueller \& Tribe \cite{dal} who introduced an alternative stochastic representation for a class of {\it linear} Cauchy problems $Lu=Fu$, with a potential $F$, including the linear wave equation. Their starting point is the well-known representation
 \be\label{DMT}
 u(t,x)
 &=&
 w(t,x)
 + \int_0^t \!\!\int_{\R^d} V(t-r,x-y)u(t-r,x-y)S(r,dy)dr,
 \ee
where $w$ is the solution with zero potential, and $S$ is a fundamental solution, restricted to be representable by a signed measure with $\sup_t |S(t,\R^d)|<\infty$. By substituting formally $u(t-r,x-y)$ on the right hand side by the last expression of $u$, one formally obtains a candidate representation for $u$ as $\sum_{m\ge 0} H_m(t,x)$, where $H_0=w$, and $H_{m}=\int_0^t \!\!\int_{\R^d} V(t-r,x-y)H_{m-1}(t-r,x-y)S(r,dy)dr$, $m\ge 1$. Finally, by convenient normalization of the kernel $S$,  the last expression induces a probabilistic representation.

\no Subsequently, Bakhtin and Mueller \cite{bak} considered the one-dimensional nonlinear wave equation
 \b*
 \partial^2_t u -\partial^2_x u
 =
 \sum_{j\ge 0} a_j u^j
 &\mbox{on}~~\R_+\times\R,&
 u(0,.)=f_1,~~\partial_tu(0,.)=f_2
 ~~\mbox{on}~~\R,
 \e*
and obtained a probabilistic representation by means of stochastic cascades, which mimics exactly the McKean  branching diffusion representation  of the KPP equation \cite{McKean}, similar to \cite{phl1,phl2}. Our  main contribution in this paper is to show that such a representation holds for a wider class of Cauchy problems, in arbitrary dimension, and with analytic nonlinearity in $(u,\partial_x u)$.

\no Our starting point is that the representations of \cite{dal} and \cite{bak} are closely related to the McKean \cite{McKean} representation of KPP equations, and the corresponding extensions in \cite{phl1,phl2,phl4}. This in fact opens the door to a much wider extension reported in Section \ref{sect:IVP} in the context of linear Cauchy problems with constant coefficients, and in Section \ref{sec:pol-u} in the context of nonlinear Cauchy problems with constant coefficients, and analytic nonlinearity in the value function $u$. The crucial tool for our extension is the so-called Duhamel formula which expresses the solution of such an equation as the convolution (i.e., integration) of the boundary conditions with respect to a family of fundamental solutions. This is in contrast with \eqref{DMT} which uses the single fundamental solution $S$, and requires that the solution $w$ of the zero-potential equation be given. Then, the probabilistic representation appears naturally after convenient normalization of the fundamental solutions.

\no It is also remarkable that the Malliavin automatic differentiation technique, exploited in \cite{phl3,phl4} in order to address semilinear parabolic second order PDEs, extends naturally to the context of general Cauchy problems by introducing the space gradient of the fundamental solution. This observation is the key-ingredient for our extension in Section \ref{sec:pol-uDu} to a general class of semilinear Cauchy problems with analytic nonlinearity in $(u,\partial_x u)$.

\no The performance of the Monte Carlo numerical method induced by our representation is illustrated in Section \ref{sect:numerics} on some relevant examples from mathematical physics. We start with two examples of semilinear wave equations: the Klein-Gordon equation in dimensions 1,2,3, which has a power nonlinearity in $u$, and a simplified version of the Yang-Mills equation (in dimension $1$), which contains a nonlinearity in the space gradient. We observe that due to some restricting conditions which will be detailed in Sections \ref{sect:IVP}, \ref{sec:pol-u}, and \ref{sec:pol-uDu}, our Monte Carlo approximation method does not apply to the multi-dimensional wave equation with space gradient nonlinearity. We next report some numerical experiments in the context of the nonlinear one-dimensional Beam equation which contains two derivatives in time and four space derivatives. We finally consider the Gross-Pitaevskii equation in dimensions $1$, $2$, and $3$, as an example of nonlinear Schr\"{o}dinger equations. All of the numerical results reveal an excellent performance of our  Monte Carlo approximation method.

\no We finally emphasize that, throughout the paper, we consider the Cauchy problem on $\R_+\times\R^d$, thus ignoring the important case of restricted domain ${\cal D}\subset\R^d$ for the space variable. We mention that Chatterjee \cite{cha} proved that the function $u(t,x)=\E_{x}[ f(tX+\sqrt{\tau}Z,B_\tau)]$, with independent r.v. $X,Z$ with standard Cauchy and normal distributions, respectively, and $\tau:=\inf\{ t>0 \;:\; B_t \notin {\cal D}\}$ is the first exit time of an independent Brownian motion from the domain $\cal D$, solves the wave equation on $\R_+\times{\cal D}$. However, this does not provide a representation for the wave equation on a restricted space domain as the determination of $f$ from given $f_1(x):=u(0, x)$ and $f_2(x):=\partial_t u(0, x)$ is not transparent.

\section{Probabilistic representation for linear Cauchy problems}
\label{sect:IVP}

\subsection{Non-homogeneous Cauchy problem}

For a smooth function $\phi:\R_+\times\R^d\longrightarrow\R$, we denote $\partial^0_t\phi=D^0\phi=\phi$, and
 \b*
 \partial^j_t \phi := \frac{\partial^j\phi}{\partial t^j}
 ~~\mbox{and}~~
 D^\alpha\phi :=\frac{\partial^{|\alpha|} \phi}{\partial x_1^{\alpha_1}\ldots \partial x_d^{\alpha_d}},
 &\mbox{for all}&
 j\ge 1,~\alpha\in\N^d.
 \e*
Given two integers $N,M\ge 1$, we denote $\N^d_M:=\{\alpha\in\N^d:|\alpha|\le M\}$, and we consider some scalar parameters $(a_j)_{1\le j\le N}, (b_\alpha)_{\alpha\in\N^d_M}\subset\R$ with
 \b*
 a_N=1
 &\mbox{and}&
 \big\{\alpha\in\N^d_M:|\alpha|=M~\mbox{and}~b_\alpha\neq 0\big\}\neq\emptyset.
 \e*
Throughout this paper, we consider nonlinear partial differential equations defined by means of the following linear Cauchy problem:
 \be
 &&
 \sum_{n=1}^N a_n \partial^n_t u
 - \sum_{\alpha\in\N^d_M} b_\alpha D^\alpha u
- F
 \;=\; 0 ~~\mbox{on}~\R_+\times\R^d,
 \label{IVP}
 \\ &&
 \partial^{n-1}_t u(0,.)=p_n f_n ~~\mbox{on}~\R^d,~~n=1,\ldots,N,
 \label{BC}
 \ee
where the boundary data and the source term satisfy the following conditions:
 \begin{equation}\label{hyp:fnF}
 \begin{array}{c}
 f_n:\R^d\longrightarrow\R,~1 \le n\le N,
 ~\mbox{and}~
 F:\R_+\times\R^d\longrightarrow\R,
 ~\mbox{are bounded continuous,}
 \end{array}
 \end{equation}
and $p_1,\ldots,p_N$ are scalar parameters in the simplex:
 \b*
 p_n>0
 &\mbox{for all}~~n=1,\ldots,N,~~\mbox{and}&
 p_1+\ldots + p_N=1.
 \e*

\subsection{Duhamel's formula for linear Cauchy problems}

We next recall the general solution of non-homogeneous Cauchy problems with constant coefficients. We first introduce the $\C^N-$valued function with components
$\hat{g}:=(\hat{g}_1,\cdots, \hat{g}_N)$:
\be\label{hatg}
 \hat g(t,\xi)
 :=
 (2 \pi)^{-{d \over 2}}e^{t B(\xi)^{\rm T}} {\rm e_1},
 &t\ge 0,~\xi\in\R^d,&
 \ee
where $( {\rm e}_1,\ldots,{\rm e}_N)$ is the canonical basis of $\R^N$,
 \begin{equation}\label{Bb}
 B(\xi)
 :=
 \left(\begin{array}{c|ccc}
        0         &          &            &              \\
        \vdots &          & I_{N-1} &              \\
        0         &         &             &              \\ \hline
        b(\xi)   & -a_1 &  \cdots & -a_{N-1} \\
        \end{array}
 \right),
 ~\mbox{and}~
 b(\xi)
 :=
 \sum_{\alpha\in\N^d_M} i^{|\alpha|}b_\alpha\xi^\alpha,
 ~
 \xi\in\R^d,
 \end{equation}
with $\xi^\alpha:=\xi_1^{\alpha_1}\ldots\xi_d^{\alpha_d}$, for all multi-index $\alpha\in\N^d$. The polynomial function $b:\R^d\longrightarrow\R$ is called {\it the symbol} of the Cauchy problem \eqref{IVP}.

\no As standard, we denote by $\Sc$ the Schwartz space of rapidly decreasing functions on $\R^d$, and by $\Sc'$ the corresponding dual space of tempered distributions. We recall that this space contains all (distributions represented by) polynomially growing functions.

\begin{Assumption}\label{hyp:hatg}
There exists $T\in(0,\infty]$ such that $\hat g(t,.)\in \Sc'$ for all $t\in [0,T)$.
\end{Assumption}


\no Under this assumption, we may introduce the so-called Green functions as the inverse Fourier transform with respect to the space variable:
 \be\label{Green}
 g(t,\cdot)
 \;:=\;
 \mathfrak{F}^{-1}\hat g(t,.),
 &t\in[0,T),&
  \ee
 in the distribution sense, i.e., $\langle \mathfrak{F}^{-1}\hat g(t,\cdot), \varphi\rangle=\langle \hat g(t,\cdot),\mathfrak{F}^{-1}\varphi\rangle$ for all $\varphi\in\Sc$, $t\in[0,T)$, where
 \b*
 \mathfrak{F}^{-1} \varphi(x)
 :=
 (2\pi)^{-d/2}\int_{\R^d}e^{i\xi\cdot x}\varphi(\xi)d\xi,
 &\mbox{for all}&
 \varphi\in\Sc.
 \e*

\begin{Assumption}\label{hyp:g}
For all $n=1,\ldots,N$:
\\
{\rm (i)} $(t,x)\longmapsto\big(g_n(t,\cdot)*\phi\big)(x)$ is continuous on $[0,T)\times\R^d$, for all bounded continuous function $\phi$ on $\R^d$;
\\
{\rm (ii)}
$g_n(t,\cdot)$ may be represented by a signed measure $g_n(t,dx)=g_n^+(t,dx)-g_n^-(t,dx)$, with total variation measure $|g_n|:=g_n^++g_n^-$ absolutely continuous with respect to some probability measure $\mu_n$; the corresponding densities $\gamma_n$, $\gamma_n^+$ and $\gamma_n^-$, defined by
 \b*
 g_n^+(t,dx)
 =
 \gamma^+_n(t,x)\mu_n(t,dx),
 ~~
 g_n^-(t,dx)
 =
 \gamma^-_n(t,x)\mu_n(t,dx),
 ~\gamma_n:=\gamma_n^+-\gamma_n^-,
 \e*
satisfy $\big|\gamma_n(t,.)\big|_\infty<\infty$, $\gamma_n(t,\R^d)<\infty$, and $\gamma_N(.,\R^d\big)\in\L^1([0,t])$ for all $t\in[0,T)$.
\end{Assumption}

\no The choice of the dominating measures $(\mu_n)_{1\le n\le N}$ will be discussed in Example \ref{eg:mu}.

\begin{Proposition}\label{prop:linearIVP}
Let $F,(f_n)_{1 \leq n \leq N}$ be as in \eqref{hyp:fnF}. Then, under Assumptions \ref{hyp:hatg} and \ref{hyp:g}, the Cauchy problem \eqref{IVP}-\eqref{BC} has a unique solution in $C^0_b([0,T]\times\R^d,\R)$ given by
 \begin{equation}\label{reprsentationconvol}
 u(t,x)
 :=
 \sum_{n=1}^{N} p_n \big(f_n*g_n(t,.)\big)(x)
+ \!\!\int_0^t \!\!(F(t-s,\cdot)* g_N(s,\cdot))(x)ds,
 ~t\in [0,T],~x\in\R^d.~
 \end{equation}
\end{Proposition}

\no We observe that the standard statement of the last proposition involves different assumptions on $g_n$, $f_n$ and $F$. Namely, one may typically relax Assumption \ref{hyp:g} and assume that $f$ and $F$ have bounded support so as to guarantee that the convolutions involved in the representation (\ref{reprsentationconvol}) are well-defined and satisfy the property $\mathfrak{F}\big(f_n*g_n(t,\cdot)\big)=(2\pi)^{d \over 2}\mathfrak{F}(f_n)\mathfrak{F}\big(g_n(t,\cdot)\big)$ for all $t\in[0,T)$, $n=1,\ldots,N$. Our conditions in Proposition \ref{prop:linearIVP} are suitable for the subsequent use in the paper.

\no For the convenience of the reader, we report the proof of Proposition \ref{prop:linearIVP}.

\proof
First, the conditions on the densities $\gamma_n$ contained in Assumption \ref{hyp:g} guarantee that the function $u$ defined in \eqref{reprsentationconvol} is bounded and continuous. Then the distribution represented by $u$ is in $\Sc'$, and we may define the corresponding Fourier transform in the space variable $\hat u(t,.):=\mathfrak{F}\big(u(t,.)\big)$ in the distribution sense. By standard calculation using the properties of the Fourier transform, we see that
 \be\label{uhat}
 \hat u(t,\xi)
 \;=\;
 \sum_{n=1}^{N} p_n\hat f_n(\xi)\,\hat g_n(t,\xi)
 +\int_0^t \hat F(t-s,\xi)\,\hat g_N(s,\xi) ds,
 &t\ge 0,~\xi\in\R^d,&
 \ee
where $\hat f:=\mathfrak{F}(f)$ and $\hat F:=\mathfrak{F}\big(F(t,.)\big)$. By the definition of $\hat g$ in \eqref{hatg}, we see that, for every fixed $\xi\in\R^d$, the function $\hat u(\cdot,\xi)$ is the unique solution of the ODE
 \be
 &&
 \sum_{n=1}^N a_n \partial^n_t \hat u
 - b(\xi)\hat u
- \hat F(.,\xi)
 \;=\;
 0 ~~\mbox{on}~\R_+,
 \label{IVP-Fourier}
 \\ &&
 \partial^{n-1}_t \hat u(0,\xi)= p_n\hat f_n(\xi),~~n=1,\ldots,N,
 \label{BC-Fourier}
 \ee
which can be written equivalently as a first order ODE in terms of the function $\hat v:=\big(\hat u,\partial_t\hat u,\ldots,\partial^{N-1}_t\hat u\big)^\top$:
 \b*
 \partial_t \hat v
 =
 B \hat v +\hat F\, {\rm e}_N,~~\mbox{on}~\R_+,
 &\mbox{and}&
 \hat{v}(0,.)=\sum_{n=1}^{N} p_n\hat f_n\,{\rm e_n}.
 \e*
Obviously, the last ODE has a unique solution with closed form obtained by the variation of the constant method $\hat{v}(t,\xi):=e^{tB(\xi)}\hat{v}(0,\xi)+\int_0^t e^{sB(\xi)}\hat F(t-s)e_N ds$, whose first entry induces the solution $\hat u$ introduced in \eqref{uhat}. To conclude the proof, it suffices to observe that any solution $\tilde u\in C^0_b([0,T)\times\R^d)$ of the Cauchy problem \eqref{IVP}-\eqref{BC} has a well-defined Fourier transform in the distribution sense satisfying the ODE \eqref{IVP-Fourier}-\eqref{BC-Fourier} for all fixed $\xi\in\R^d$.
\ep


\begin{Example}[Heat equation] Let $N=1$, $M=2$, $b_\alpha=0$ whenever $|\alpha|\le 1$, and $b_{1,1}=b_{2,2}=1$, $b_{1,2}=b_{2,1}=0$. Then $B(\xi)=b(\xi)=-|\xi|^2,$ $\xi\in\R^d$, and
 \b*
 g_1(r,z)
 := (2\pi)^{-{d}}\int e^{-|\xi|^2 r +i\xi z}d\xi
 = (4\pi r)^{-d/2}e^{-\frac{|z|^2}{4r}},
 &(r,z)\in\R_+\times\R^d.&
 \e*
\end{Example}

\begin{Example}[Airy equation] Let $N=1$, $M=3$, $b_0=0$, $b_\alpha=0$ whenever $|\alpha|\le 2$, and $b_\alpha=0$, and $\xi\longmapsto b(\xi)$ odd in the sense that
 \b*
 &
  b\big(-\xi_j{\rm e}_j+\sum_{\ell\neq j} \xi_\ell{\rm e}_\ell\big)
  =
 -b_\alpha\big(\xi_j{\rm e}_j+\sum_{\ell\neq j} \xi_\ell{\rm e}_\ell\big)
 &
 \mbox{for all}~~j=1,\ldots,d.
 \e*
Then $B(\xi)=b(\xi)$ is scalar valued, and
 \b*
 g_1(r,z)
 :=  (2\pi)^{-{d}}\int e^{i(-r b(\xi) + \xi z)}d\xi
 = \frac{1}{(3r)^{d/3}} {\rm Ai}_b\Big(\frac{z}{(3r)^{1/3}}\Big),
 &(r,z)\in\R_+\times\R^d,&
 \e*
where we introduced the $d-$dimensional Airy function
 \b*
 {\rm Ai}_b(x) :=  (\pi)^{-{d}} \int_{\R_+^d} \cos{\big( -{b(\xi) \over 3}+x\cdot\xi\big)} d\xi,
 &x\in\R^d.&
 \e*
\no For $d=1$, notice that $\int |g_1|(t,dz)=\infty$, so that Assumption \ref{hyp:g} fails in this example.
\end{Example}

\begin{Example}[Wave equation] Let $N=2$, with $a_1=0$, $M=2$, $b_\alpha=0$ whenever $|\alpha|\le 1$, and $b_{1,1}=b_{2,2}=1$, $b_{1,2}=b_{2,1}=0$. Then
 $$
 B(\xi)=\left(\begin{array}{cc}
                  0         & 1 \\
                  -|\xi|^2 & 0
                  \end{array}
            \right),
 ~~
 e^{B(\xi)r}=(2\pi)^{-{d \over 2}}\left(\begin{array}{cc}
                                     \cos{(r|\xi|)}       & |\xi|^{-1}\sin{(r|\xi|)} \\
                                     -|\xi|\sin{(r|\xi|)} & \cos{(r|\xi|)}
                  \end{array}
            \right),
 ~~\xi\in\R^d,
 $$
and, for $(r,z)\in\R_+\times\R^d$,
 $$
 g_2(r,dz)
 =
 (2\pi)^{-d}\int \frac{\sin{(r|\xi|)}}{|\xi|}e^{i\xi\cdot z}d\xi,
 ~
 g_1(r,dz)
 =
 (2\pi)^{-d}\int \cos{(r|\xi|)}e^{i\xi\cdot z}d\xi
 =
 \partial_r g_2(r,dz).
 $$
We have (see Kirchhoff's formula in \cite{eva} for closed-expression  for $g_1$ and $g_2$ in $\RR^d$)
 \beaa
 g_2(r,dz)
 &=&
 \left\{\begin{array}{ll}
         {1 \over 2}\,1_{\{|z|<r\}}dz, & \mathrm{for}~ d=1
         \\
         {1 \over 2 \pi} (r^2-|z|^2)^{-\frac12} \,1_{\{|z|<r\}}dz, & \mathrm{for}~ d=2
         \\
         {\sigma_r (dz) \over 4 \pi r}, & \mathrm{for} \; d=3
         \end{array}
         \right.
 \eeaa
where $\sigma_r(dz)$ denotes the surface area on $\partial B(0,r)$. Assumption \ref{hyp:g} are discussed in Example \ref{eg:mu}.
\end{Example}

\begin{Example}[Beam equation] Let $d=1$, $N=2$, with $a_1=0$, $M=4$, $b_1=b_2=b_3=0$, and $b_4=1$ corresponding to the fourth-order PDE $u_{tt} +\partial_{x}^4 u=0$. Then,
 \beaa
 g_1(r,z)=\partial_r g_2(r,z)
 &\mbox{and}&
 g_2(r,z)=\sqrt{r} G\left({z \over \sqrt{r}}\right),
 \eeaa
where $G(0)=1/\sqrt{2 \pi}$, and
 \be\label{Gbeam}
 2G'(x)
 &=&
 \mathrm{F}_s\left(\frac{x}{\sqrt{2 \pi
   }}\right)-\mathrm{F}_c\left(\frac{x}{\sqrt{2 \pi }}\right), \quad x \in \RR,
 \ee
with the Fresnel integrals $\mathrm{F}_c(x):=\int_0^x \cos(\pi t^2/2) dt$ and $\mathrm{F}_s(x):=\int_0^x \sin(\pi t^2/2) dt$. Note that $\int_{\RR} |G(x)| dx <\infty$ as $G(x) \underset{|x| \rightarrow \infty}{\sim} \frac{\sqrt{\frac{2}{\pi }} \left(\cos \left(\frac{x^2}{4}\right)-\sin
   \left(\frac{x^2}{4}\right)\right)}{x^2}$.
\end{Example}

\begin{Remark}\label{rem:CharacteristicPolynomial}
{\rm In order to compute the Green functions $g$ as the inverse Fourier transform of the associated $\hat g$, one needs to diagonalize the matrix $B(\xi)$. Direct examination reveals that the eigenvalues of the matrix $B(\xi)$ are the roots of the corresponding characteristic polynomial $b(\xi)-\sum_{n=1}^{N}a_n\lambda^n $. Assume that $B(\xi)$ has $N$ distinct (simple) eigenvalues $\big(\lambda_j(\xi)\big)_{1\le n\le N}\in\C^N$. Then, denoting by diag$[\lambda]$ the diagonal matrix with diag$[\lambda]_{j,j}=\lambda_j$, it follows that
 \b*
 B(\xi)
 =
 P(\xi){\rm diag}[\lambda(\xi)]P(\xi)^{-1},
 &\mbox{where}&
 P(\xi)_{j,\ell} := \lambda_\ell(\xi)^j,~~1\le j,\ell\le N.
 \e*
The matrix $P(\xi)$ is the so-called Vandermonde matrix whose inverse is given by:
 \b*
 \big\{P(\xi)^{-1}\big\}_{j,n}
 =
 \frac{\Lambda_n(\xi)}{\lambda_j(\xi)\prod_{\ell\neq j}(\lambda_\ell-\lambda_j)(\xi)},
 &j,n=1,\ldots,N,&
 \e*
where
 \b*
 \Lambda_N(\xi)=1,
 &\mbox{and}&
 \Lambda_n(\xi)
 :=
 (-1)^{n-1}\!\!\!\!\sum_{_{\!\!\!\tiny\begin{array}{c}
                    1\!\!\le\!\!\ell_1\!\!\le\!\!\ldots\!\!\le\!\!\ell_{N-n}\!\!\le\!\! N
                                                                                      \\
                                                                                      \ell_1,\ldots,\ell_{N-n}\!\!\neq\!\! n
                    \end{array}}}
                    \!\!\!\!(\lambda_{\ell_1}\cdots\lambda_{\ell_{N-n}})(\xi)
 ~~\mbox{for}~n<N.
 \e*
Therefore, by the definition of $\hat g_n$ in \eqref{hatg}, we have for $(r,\xi)\in\R_+\times\R^d$:
 \be\label{hatgn:rem}
 \hat g_n(r,\xi)
 &=&
 (2\pi)^{-{d \over 2}}
 \Lambda_n(\xi)\sum_{j=1}^N \frac{e^{r\lambda_j(\xi)}}
                                                       {\prod_{\ell\neq j}(\lambda_\ell-\lambda_j)(\xi)},
 ~n=1,\ldots,N.
 \ee
\ep}
\end{Remark}

\begin{Remark}{\rm
All of the results of the present paper extend to the case of Cauchy problems for complex-valued functions, with coefficients $(a_n)_{1\le n<N}$ and $(b_\alpha)_{\alpha\in\N^d_M}$ in $\C$, thus allowing to include, for instance, the Schr\"odinger equation (see Example \ref{Schrodinger equation and analytical continuation}) and its semilinear extension as the Gross--Pitaevskii PDE (see Section \ref{Gross--Pitaevskii}). This extension follows by simply applying the methodology described throughout the paper separately to the real part and the imaginary part of the representation.
\ep}
\end{Remark}

\subsection{Probabilistic representation}

Let $(\Omega,\Fc,\P)$ be a probability space supporting two random variables $\tau$ and $I$, with
 \begin{equation}\label{Itau}
 \tau~\mbox{and}~I~\mbox{independent,}
 ~\P[\tau\in dt]=\rho(t)\1_{\{t\ge 0\}} dt
 ~\mbox{and}
 ~\P[I=n]=p_n,
 ~n=1,\ldots,N,
 \end{equation}
for some $C^0(\R_+,\R)$ density function $\rho>0$ on $(0,\infty)$. We shall denote $\bar\rho(t):=\int_t^\infty \rho(s)ds$.

\no Recall the densities $\gamma_n$ and the dominating probability measures $\mu_n(t,\cdot)$, $n=1,\ldots,N$, introduced in Assumption \ref{hyp:g}. For all $t\ge 0$, we introduce the random variables
 \begin{equation}
 X^n_t:=x+Z^n_t
 ~\mbox{independent of}~(I,\tau),~\mbox{with}~\P\big[Z^n_t\in dz]=\mu_n(t,dz),
 ~n=1,\ldots,N.
 \end{equation}

\begin{Example}\label{eg:mu}
As $g_n(t,\cdot)\in\L^1(\R^d)$ by Assumption \ref{hyp:g}, we may choose $\mu_n(t,dz)=\|g_n(t,\cdot)\|_{\L^1(\R^d)}^{-1} |g_n|(t,dz)$ and $\gamma_n(t,z)=\mbox{\rm sgn}(g_n)(t,z) \|g_n(t,\cdot)\|_{\L^1(\R^d)}$, where $\mbox{\rm sgn}(\alpha):=\1_{\{\alpha\ge 0\}}-\1_{\{\alpha<0\}}$.

\no {\rm(i)} For the heat equation, we directly compute that $\|g_1(t,\cdot)\|_{\L^1(\R^d)}=1$. Then, choosing $\gamma_1(t,\cdot)= 1$ leads to $ Z_t^1=\sqrt{2 t} Z$ with $Z \in N(0,I_d)$ a $d$-Gaussian random variable.

\no {\rm(ii)} For the wave equation, we directly compute that $\|g_2(t,\cdot)\|_{\L^1(\R^d)}=t$. Then, choosing $\gamma_2(t,\cdot)= t$ leads to $ Z_t^2=t Z$ where

$\bullet$ in dimension $d=1$, $Z$ has a uniform distribution on $[-1,1]$

$\bullet$ in dimension $d=2$, the law of $Z$ is defined by the density ${1 \over 2 \pi}{1 \over \sqrt{1-z^2}}1_{|z|<1}$,

$\bullet$ in dimension $d=3$, the law of $Z$ is ${1 \over 4\pi} \mu_{S^2}(dz)$, where $\mu_{S^2}$ denotes the volume measure on the unit sphere.

\no \rm{(iii)} For the Beam equation, $\gamma_2(t,z)=t\; \mathrm{sgn}(G(z)) \|G\|_{\L^1(\R^d)}$ where $G$ is defined in \eqref{Gbeam}. We have $Z_t^2= \sqrt{t} Z $ with $Z$ distributed according to $\|G\|_{\L^1(\R^d)}^{-1}|G(z)|$.
\end{Example}

\begin{Example}[Schr\"odinger equation and analytical continuation]The Duhamel formula for the Schr\"odinger equation of a free particle (with source term)
\bea i \partial_t u =-{1 \over 2}\Delta u +F, \quad u(0,x)=f_1(x) \label{DuhamelQuantum} \eea is
\beaa u(t,x)= \int_{\RR^d} dy f_1(y) g_1(t,x-y)-i
\int_0^t \int_{\RR^d} dy g_1(s,x-y) F(s,y) ds \eeaa with
$ g_1(t,x)={e^{- {x^2 \over 2 i t}} \over (2 \pi i t)^{d \over 2} }$. Note the coefficient $-i$ in front of $F$ as $a_1=i$ here. By setting $y-x=\sqrt{i t}z$, this can be written as
\beaa u(t,x)=  \int_{-\infty e^{-{i \pi \over 4}}}^{\infty e^{-{i \pi \over 4}}} dz f_1(x+e^{{i \pi \over 4}} \sqrt{ t}z) e^{-{z^2 \over 2}}  -i \int_0^t ds \int_{-\infty e^{-i \pi \over 4}}^{\infty e^{-i \pi \over 4}}
 F(s,x+e^{{i \pi \over 4}} \sqrt{s}z)e^{-{z^2 \over 2}} dz  ds  \eeaa By assuming that $|f_1(x+e^{{i \pi \over 4}+i\theta} \sqrt{ t}R) |e^{-{R^2 \over 2}\cos (2\theta)}$ and $|F(s,x+e^{{i \pi \over 4}+i\theta} \sqrt{ s}R) |e^{-{R^2 \over 2}\cos (2\theta)}$ goes to $zero$ as $R \rightarrow \infty$ when $\theta \in [-{\pi \over 4},0]^d \cup [{3\pi \over 4},\pi]^d $, the integration over $z$ in $[-\infty e^{-{i \pi \over 4}},\infty e^{-{i \pi \over 4}}]^d$ can be deformed into
 $[-\infty ,+\infty ]^d$ by analytical continuation  and we obtain
 \beaa u(t,x)=  \EE[ f_1(x+e^{{i \pi \over 4}} \sqrt{ t}Z) ]  -i
\int_0^t ds \EE[ F(s,x+e^{{i \pi \over 4}} \sqrt{s}Z)]  ds  \eeaa with $Z \in N(0,I_d)$. \label{Schrodinger equation and analytical continuation}
\end{Example}

\no The following representation result is a simple rewriting of Proposition \ref{prop:linearIVP} in terms of the last notations.

\begin{Proposition}
Let $(f_n)_{1 \leq n \leq N}$ and $F$ be as in \eqref{hyp:fnF}. Then, under Assumptions \ref{hyp:hatg} and \ref{hyp:g}, the unique $C^0_b([0,T]\times\R^d,\R)$ solution of the Cauchy problem \eqref{IVP}-\eqref{BC} is:
 $$
 u(t,x)
 =
 \E\Big[ \1_{\{\tau\ge t\}}\, \frac{\gamma_{I}(t,Z^I_t)}{\bar\rho(t)} f_I\big(X^I_\tau\big)
            + \1_{\{\tau<t\}}\, \frac{\gamma_N(\tau,Z^N_\tau)}{\rho(\tau)} F\big(t-\tau, X^N_\tau\big)
     \Big],
 ~t<T,x\in\R^d.
 $$
\end{Proposition}

\section{A first class of semilinear Cauchy problems}
\label{sec:pol-u}

In this section, we consider the semilinear Cauchy problem:
\be
 &&
 \sum_{n=1}^N a_n \partial^n_t u
 - \sum_{\alpha\in\N^d_M} b_\alpha D^\alpha u
 - \sum_{j\ge 0} q_j\, c_j\, u^j
 \;=\;
 0 ~~\mbox{on}~\R_+\times\R^d,
 \label{IVP1}
 \\ &&
 \partial^{n-1}_t u(0,.)=p_nf_n ~~\mbox{on}~\R^d,~~n=1,\ldots,N,
 \label{BC1}
 \ee
where the nonlinearity is defined by means of an atomic probability measure $(q_j)_{j\ge 0}$, with $q_j\ge 0$ and $\sum_{j\ge 0}q_j=1$, together with the functions
 \be\label{cond:cell}
 c_j:\R_+\times\R^d \longrightarrow\R,
 &j\ge 0,&
 \mbox{bounded continuous, for all}~j\ge 0.
 \ee
In order for the power series on the right hand-side of \eqref{IVP1} to be well defined, we assume that
 \be\label{cond:F}
 H_1(s)
 &:=&
 \sum_{j\ge 0} \,q_j\|c_j\|_\infty\, s^j.
 \ee
has a strictly positive radius of convergence, so that it is well defined at least in some neighborhood of the origin.

\subsection{The branching mechanism}
\label{sect:branching}

A particle of generation $\nu\in\N$ is a multi-integer $k:=(0,k_1,\ldots, k_\nu)\in\mathbf{K}_\nu:=\{0\}\times\N^\nu$. We set by convention $\mathbf{K}_0:=\{0\}$, and we denote by $\mathbf{K}:=\cup_{\nu\ge 0} \mathbf{K}_\nu$ the collection of all particles. For $\nu\ge 1$, and a particle $k:=(0,k_1,\ldots, k_\nu) \in \mathbf{K}_\nu$, we denote by $k_-$ its parent particle
 \b*
 k_-:=(0,k_1,\ldots, k_{\nu-1}) \in \mathbf{K}_{\nu-1}.
 \e*
We next introduce independent families of random variables $(\tau_k)_{k\in\mathbf{K}}$, $(I^0_k)_{k\in\mathbf{K}}$, and $(J_k)_{k\in\mathbf{K}}$:
\begin{itemize}
\item $I^0_k$ and $\tau_k$ are iid copies of the random variables $I$ and $\tau$, respectively, as introduced in \eqref{Itau};
\item $J_k$ are iid random variables with $\P[J_k=j]=q_j$ for all $j\ge 0$.
\end{itemize}
The time occurrences of the branching events, prior to $t$, are recorded through the sequence $(T^t_k)_k$ defined for all $t\ge 0$ by:
 \b*
 T^t_{0-}:=0,
 ~~\mbox{and}~
 T^t_k
 :=
 t\wedge\big(T^t_{k_-}+\tau_k\big),
 &\mbox{for all}&
 k\in\mathbf{K}.
 \e*
With these notations, each particle $k$ lives on the time interval $[T^t_{k_-},T^t_k]$. The branching mechanism is then the following:
\begin{itemize}
\item Start from particle $0$ at the time origin.
\item At the first branching time $T^t_0=t\wedge\tau_0$, particle zero dies out; if $T^t_0<t$, it generates $J_0$ descendants labelled $(0,1),\ldots,(0,J_0)$.
\item Each descendant particle $k$ undergoes the same behavior, independently of its peer-particles: it dies out at the branching time $T^t_k$, and generates $J_k$ descendants labelled $(k,1),\ldots, (k,J_k)$, whenever $T^t_k<t$.
\item We denote by $\Kc^t_s$ the collection of all living particles at time $s$, and by $\overline{\Kc}_t:=\cup_{s\le t}\Kc^t_s$ the collection of all particles which have been living prior to time $t$. For simplicity we set $\Kc_t:=\Kc^t_t$.
\item We finally denote for all particle $k\in\overline{\Kc}_t$:
 \be\label{type-1st}
 I^t_k
 &:=&
 I^0_k \1_{k\in\Kc_t}
 + N \1_{k\in \overline{\Kc}_t\setminus\Kc_t}.
 \ee
\end{itemize}

\subsection{Probabilistic representation}

Given the branching mechanism defined in the previous subsection, we introduce the corresponding  branching process:
 \be\label{Xk}
 X^0_0=x,
 &\mbox{and}&
 X^k_{T^t_k}
 :=
 X^{k_-}_{T^t_{k_-}} + Z^k_{T^t_k},
 ~~\mbox{for all}~~
 k\in\overline{\Kc}_t,
 \ee
where the distribution of $Z^k_{T^t_k}$ depends on the type of the particle:
 \be\label{Zk}
 \P\big[ Z^k_{T^t_k}\in dz \big| I^t_k, \Delta T^t_k \big]
 \;=\;
 \mu_{I^t_k}(\Delta T^t_k, dz),
 &\mbox{with}&
 \Delta T^t_k:=T^t_k-T^t_{k_-}.
 \ee
We now introduce, for all $t\ge 0,$ and $x\in\R^d$, the random variable
 \begin{equation}\label{xi-semilinear1}
 \xi_{t,x}
 :=
 \prod_{k\in\Kc_t} \frac{\gamma_{I^t_k}\big(\Delta T^t_k,Z^k_{T^t_k}\big)}
                                    {\bar\rho\big(\Delta T^t_k\big)}
                            f_{I^t_k}\big(X^{k}_t\big)
 \prod_{k\in\overline{\Kc}_t\setminus \Kc_t}\frac{\gamma_N\big(\Delta T^t_k,Z^k_{T^t_k}\big)}
                                    {\rho\big(\Delta T^t_k\big)}
                                                                    c_{J_k}\big(t-T^t_k,X^k_{T^t_k}\big).
 \end{equation}



\no Recall the power series $H_1$ defined in \eqref{cond:F}.

\begin{Assumption}\label{assum:semilinear1}
The power series $H_1$ has a radius of convergence $R_1\in(0,\infty]$, and

{\rm (i)} $r_1:=\sup_{1\le n\le N}
                    \|f_n\|_\infty \|\gamma_n\|_\infty < R_1$;

{\rm (ii)} there are constants $T>0$ and $s_1>r_1$ such that $\int_{r_1}^{s_1} H_1(s)^{-1}ds=T\|\gamma_N\|_\infty$.
\end{Assumption}

\begin{Theorem}\label{thm:IVP1}
Let $f\in C^0_b(\R^d,\R^n)$ and $c_j\in C^0_b(\R_+\times\R^d,\R)$. Then, under Assumptions \ref{hyp:hatg}, \ref{hyp:g} and \ref{assum:semilinear1}, we have $\xi_{t,x}\in\L^1$ and $u(t,x):=\E\big[\xi_{t,x}\big]$ is the  unique solution of the semilinear Cauchy problem \eqref{IVP1}-\eqref{BC1}.
\end{Theorem}

\proof {\bf 1.} We first verify the $\L^1-$integrability of $\xi_{t,x}$. From the expression of $\xi_{t,x}$ in \eqref{xi-semilinear1}, we see that
 \be\label{chi-1}
 \big|\xi_{t,x}\big|
 &\le&
 \chi_t
 \;:=\;
 \prod_{k\in\Kc_t} \frac{r_1}{\bar\rho\big(\Delta T^t_k\big)}
 \prod_{k\in\overline{\Kc}_t\setminus\Kc_t} \frac{\overline{c}_{J_k}\|\gamma_N\|_\infty}
                                                                     {\rho\big(\Delta T^t_k\big)}.
 \ee
Our objective is now to prove that $t\longmapsto\E \chi_t<\infty$. By standard arguments, this function is related to the ODE
 \be\label{ODE0}
 \partial_t w \;=\; \|\gamma_N\|_\infty H_1(w)
 &\mbox{and}&
 w(0)=r_1.
 \ee
Let us first verify that this ODE has a non-exploding solution under the second condition in Assumption \ref{assum:semilinear1}. Since the power series function $H_1$ has radius $R_1$, the condition $r_1<R_1$ is necessary to find a (finite) solution of the last ODE. Next, for an arbitrary $L>0$, notice that $H_1$ is Lipschitz on $[-L,L]$, so that the last ODE has a unique non-negative solution $w$ up to $T_L$, where
 \b*
 \lim_{t\to T_L} w(t) = L
 &\mbox{whenever}&
 T_L<\infty,
 \e*
and the non-negativity of $w$ follows from the fact that $w(0)>0$ and all derivatives $H^{(i)}_1(0)\ge 0$ for all $i\ge 0$. Then, it follows from direct integration of the ODE that
 \b*
 \|\gamma_N\|_\infty T_L
 =
 \int_0^{T_L} \frac{\partial_s w}{H_1(w)}(s)ds
 =
 \int_{r_1}^{L} \frac{ds}{H_1(s)}
 &\mbox{as long as}&
 T_L<\infty.
 \e*
This proves that $w$ is non-exploding solution on $[0,T]$ if and only if $\int_{r_1}^{s_1} \frac{ds}{H_1(s)}=\|\gamma_N\|_\infty T$ for some constant $s_1>r_1$.

Under the last condition, we obtain by direct integration of the ODE \eqref{ODE0} that:
 \b*
 w(t)
 &=&
 r_1+\|\gamma_N\|_\infty\sum_{j\ge 0}q_j\bar c_j \int_0^t w(s)^j ds
 \\
 &=&
 \E \Big[ \frac{r_1\1_{\{\tau_0\ge t\}}}{\bar\rho(t)}
              + \1_{\{\tau_0<t\}} \frac{\|\gamma_N\|_\infty\|c_{J_0}\|_\infty}{\rho(\tau_0)} w(\tau_0)^{J_0}
     \Big]
 \\
 &=&
 \E \Big[ \frac{r_1\1_{\{\tau_0\ge t\}}}{\bar\rho(t)}
              + \1_{\{\tau_0<t\}} \frac{\|\gamma_N\|_\infty\|c_{J_0}\|_\infty}{\rho(\tau_0)}
 \\ && \hspace{37mm}
                                           \prod_{k\in\Kc^1_t} \Big(\frac{r_1\1_{\{T^t_k\ge t\}}}{\bar\rho(\tau_k)}
                                                                                  + \1_{\{T^t_k<t\}} \frac{\|\gamma_N\|_\infty\|c_{J_0}\|_\infty}
                                                                                                                     {\rho(\tau_k)}
                                                                                                             w(T^t_k)^{J_k}
                                                                         \Big)
     \Big],
 \e*
where $\Kc^n_t$ denotes the particles generated after $n$ defaults prior to $t$ and the last two equalities follow from the definition of the branching mechanism together with the tower property. Iterating up to the $n-$th generation, this provides:
 \b*
 v(t)
 =
 \E \chi_t^n,
 &\mbox{where}&
 \chi_t^n
 :=
 \prod_{k\in\Kc^n_t}\frac{r_1}{\bar\rho(\Delta T^t_k)}
 \prod_{k\in\cup_{j<n} \Kc^j_t} \frac{\|\gamma_N\|_\infty\|c_{J_0}\|_\infty}{\rho(\Delta T^t_k)}
 \\
 && \hspace{20mm}
 + \prod_{k\in\cup_{j\le n} \Kc^j_t} \frac{\|\gamma_N\|_\infty\|c_{J_0}\|_\infty}{\rho(\Delta T^t_k)}
    \prod_{k\in\Kc^{n+1}_t} w(T^t_{k-})^{J_{k-}}.
 \e*
As $w\ge 0$, it follows that $\chi^n_t\ge 0$, and we then deduce from Fatou's lemma that
 \b*
 w(t)
 =
 \liminf_{n\to\infty} \E \chi_t^n
 \ge
 \E \chi_t.
 \e*
Since $w(t)<\infty$, this provides the required integrability of the bound $\chi_t$ in \eqref{chi-1}.

\no {\bf 2.} We next prove that $v(t,x):=\E\big[\xi_{t,x}\big]$ is the unique solution of the nonlinear Cauchy problem \eqref{IVP1}. To see this, observe that
 \b*
 \xi_{t,x}
 &=&
 \1_{\{\tau_0\ge t\}}\frac{\gamma_{I_0^t}(t,Z^0_t)}
                                     {\bar\rho(t)}
                              f_{I_0^t}(X^0_t)
 +\1_{\{\tau_0<t\}}\frac{\gamma_N(\tau_0,Z^0_{\tau_0})}
                                     {\rho(\tau_0)}
                              c_{J_0}(t-\tau_0,X^0_{\tau_0})
                             \prod_{j=1}^{J_0} \xi^{(j)}_{t-\tau_0,X^0_{\tau_0}},
 \e*
where $\xi^{(j)}_{t-\tau_0,X^0_{\tau_0}}$ have the same distribution, conditional on $\big(\tau_0, X^0_{\tau_0}\big)$. Taking expectations, and using the tower property, it follows from the definition of the function $v$ that
 \b*
 v(t,x)
 &=&
 \E\Big[ \1_{\{\tau_0\ge t\}}\frac{\gamma_{I_0^t}(t,Z^0_t)}
                    {\bar\rho(t)}
             f_{I_0^t}(X^0_t)
             +\1_{\{\tau_0<t\}}\frac{\gamma_N(\tau_0,Z^0_{\tau_0})}
                                               {\rho(\tau_0)}
                                        c_{J_0}(t-\tau_0,X^0_{\tau_0})
                                        v\big(t-\tau_0, X^0_{\tau_0}\big)^{J_0}
     \Big]
 \\
 &=&
 \sum_{n=1}^N \int p_n(f_n*g_n(t,.)\big)(x)
 + \int_0^t \big(F(t-s,\cdot)*g_N(s,\cdot)\big)(x) ds,
 \e*
where $F:=\sum_{j\ge 0} q_j c_j v^j$. Since $v$ is bounded, it follows from Proposition \ref{prop:linearIVP} that $v$ is the unique $C^0_b([0,T]\times\R^d,\R)$ solution of the nonlinear Cauchy problem \eqref{IVP1}-\eqref{BC1}.
 \ep

\begin{Remark}[Finite propagation speed] {\rm From the simulation of $Z_t^2$ in the case of the wave equation, we deduce directly the finite propagation speed property: If $f_1=f_2=0$ on a ball $B(x_0,t_0)$ of center $x_0$ and radius $t_0$ then $u=0$ within the cone $K(x_0,t_0):=\big\{ (x,t):~0 \leq t \leq t_0~\mbox{and}~|x-x_0| \leq t_0-t\big\}$.
\ep}
\end{Remark}

\no The next statement provides some sufficient conditions for $\xi_{t,x}$ to have a finite $p-$th moment. In particular, for $p=2$, this guarantees that the error estimate of the Monte Carlo approximation induced by the representation of Theorem \ref{thm:IVP1} is characterized by the standard central limit theorem, and is in particular independent of the dimension $d$ of the space variable $x$. Our first requirement is stated on the power series:
 \be\label{cond:Fp}
 H_p(s)
 &:=&
 \sum_{j\ge 0} (q_j \|c_j\|_\infty)^p s^j.
 \ee

\begin{Assumption}\label{assum:semilinear1-L2} The power series $H_p$ has a radius of convergence $R_p\in(0,\infty]$, and

{\rm (i)}  $r_p:=\max_{1\le n\le N}
          \|f_n\|_\infty^p \big\|\gamma_n\bar\rho^{\frac{1-p}{p}}\big\|_\infty^p<R_p$, and $\alpha_p:=\big\| \gamma_N\rho^{\frac{1-p}{p}} \big\|_\infty^p < \infty$;

{\rm (ii)} there are constants $T>0$ and $s_p>r_p$ such that $\int_{r_p}^{s_p}H_p(s)^{-1}ds=\alpha_p T$.
\end{Assumption}

\begin{Theorem}\label{Thmsemilinear}
Let $f\in C^0_b(\R^d,\R^n)$ and $c_j\in C^0_b(\R_+\times\R^d,\R)$. Then, under Assumptions \ref{hyp:hatg}, \ref{hyp:g} and \ref{assum:semilinear1-L2}, we have $\xi_{t,x}\in\L^p$ for all $t\in[0,T]$.
\end{Theorem}

\proof Similar to the calculation in the previous proof, we have
 \b*
 \big|\xi_{t,x}\big|^p
 &\le&
 \prod_{k\in\Kc_t} \frac{r_p}{\bar\rho(\Delta T_k^t)}
 \prod_{k\in\overline{\Kc}_t\setminus\Kc_t} \frac{\alpha_p\|c_{J_k}\|_\infty^p}{\rho(\Delta T_k^t)}.
 \e*
The right hand-side is now related to the ODE
 \b*
 \partial_t w \;=\; \alpha_p H_p(w)
 &\mbox{and}&
 w(0)=r_p.
 \e*
the integrability of $|\xi_{t,x}|^p$ can now be verified by following the same line of argument as in Step 1 of the proof of Theorem \ref{thm:IVP1}.
 \ep

\section{Further nonlinear Cauchy problems}
\label{sec:pol-uDu}

In this section, we consider the following semilinear Cauchy problem with polynomial nonlinearity in the pair $(u,Du)$:
\be\label{IVP2}
 &&
 \sum_{n=1}^N a_n \partial^n_t u
 - \sum_{\alpha\in\N^d_M} b_\alpha D^\alpha u
 - \sum_{j\ge 0} q_j\, c_{j,0}\, u^{\ell_{j,0}} \prod_{h=1}^H(c_{j,h}\cdot Du)^{\ell_{j,h}}
 \;=\;
 0~~\mbox{on}~\R_+\times\R^d
 \label{IVP2}
 \\ &&
 \partial^{n-1}_t u(0,.)=p_nf_n ~~\mbox{on}~\R^d,~~n=1,\ldots,N,
 \label{BC2}
 \ee
where $(\ell_j)_{j\ge 0}\subset\N^{1+H}$ is a sequence of vector integers $\ell_j=(\ell_{j,0},\ldots,\ell_{j,H})$,
and the nonlinearity is defined by means of an atomic probability measure $(q_j)_{j\ge 0}$, with $q_j>0$ and $\sum_{j\ge 0}q_j=1$, together with the functions
 \begin{equation}\label{cond:cell-bh}
 c_{j,0}:\R_+\times\R^d \longrightarrow\R,
 ~\mbox{bounded continuous, and}~
 c_{j,h}:\R_+\times\R^d \longrightarrow B_1(\R^d),
 j\ge 0,~1\le h\le H.
 \end{equation}
Here, $B_1(\R^d)$ denotes the unit ball in $\R^d$ in the sense of the Euclidean norm. Notice that the above Cauchy problem covers the particular case of \eqref{IVP1}-\eqref{BC1} by setting $\ell_{j,h}=0$ for all $h=1,\ldots,H$.

\no Recall the Green functions $g=(g_n)_{1\le n\le N}$ from \eqref{Green}. By standard Fourier transform theory, we have
 \b*
 \partial_xg_n(t,\cdot)
 &=&
 i\mathfrak{F}^{-1}(\xi\mapsto\xi\hat g_n(t,\xi)\big),
 ~~t\ge 0,~n=1,\ldots,N.
 \e*
Similar to the probability measures $\mu_n$ introduced in Assumption \ref{hyp:g}, we now assume that the distributions $\partial_xg_n$ can be represented by signed measures for which we may introduce dominating measures for $|\partial_xg_n(t,\cdot)|:=\sum_{m=1}^d|\partial_{x_m}g_n(t,\cdot)|$.

\begin{Assumption}\label{hyp:Dgn}
{\rm (i)} For all $n=1,\ldots,N$ and $m=1,\ldots,d,$ the distribution $\partial_{x_m}g_n(t,\cdot)$ may be represented by a signed measure with total variation $|\partial_{x_m}g_n(t,\cdot)|$.
\\
{\rm (ii)} The measure $|\partial_xg_n(t,\cdot)|$ is absolutely continuous with respect to some probability measure $\mu_n^1(t,\cdot)$ so that we may define the density vector $\gamma^1_n=\big(\gamma^1_{n,1},\ldots,\gamma^1_{n,d}\big)$ by
 \b*
 \partial_xg_n(t,dx)
 &=&
 \gamma^1_n(t,x)\mu^1_n(t,dx),
 ~~t\ge 0,~x\in\R^d.
 \e*
\end{Assumption}

\no Notice that although $\gamma^1_n(t,.)$ is defined $d\mu^1_n(t,\cdot)-$a.s. this will be sufficient for our needs. Our starting point is the following ``automatic differentiation property'', which follows by direct differentiation of the Duhamel formula of Proposition \ref{prop:linearIVP}.

\begin{Proposition}
In addition to the conditions of Proposition \ref{prop:linearIVP}, let Assumption \ref{hyp:Dgn} hold true. Then, the solution $u$ of the linear Cauchy problem \eqref{IVP}-\eqref{BC} is differentiable with respect to the space variable with
 \b*
 \partial_x u(t,x)
 &=&
 \E\Big[ \1_{\{\tau\ge t\}}\, \frac{\gamma^1_{I}(t,Z^{1,I}_t)}{\bar\rho(t)} f_I\big(X^{1,I}_\tau\big)
            + \1_{\{\tau<t\}}\, \frac{\gamma^1_N(\tau,Z^{1,N}_\tau)}{\rho(\tau)}\, F\big(t-\tau, X^{1,N}_\tau\big)
     \Big],
 \e*
where $\tau$ is a random time with density $\rho$, and $X^{1,I}_\tau=x+Z^{1,I}_\tau$, with $Z^{1,I}_t$ distributed as $\mu^1(t,\cdot)$ independent of $\tau$.
\end{Proposition}

\begin{Example}
Let us illustrate the last result on our main examples.

\no {\rm (i)} Heat equation: we directly compute that $\partial_x g_1(t,dx)=-{x \over 2t}g_1(t,dx)$.

\no {\rm (ii)} Wave equation, $d=1$: we have $\partial_x g_2(t,dx)={1 \over 2}\left( \delta(x+t)-\delta(x-t) \right)dx$, but $\partial_x g_1=\partial^2_{xx} g_2$ can not be represented as a signed measure, thus violating Assumption \ref{hyp:Dgn}. However, we may still handle the one dimensional wave equation by reducing to the case $f_1=0$, see Section \ref{sect:YangMills}.

\no {\rm (iii)} Wave equation, $d>1$: Assumption \ref{hyp:Dgn} is not satisfied as $\partial_x g_1$ and $\partial_x g_2$ involve first-order derivative of the delta function supported on the light cone $\{ (x,t) \;:\; t^2-x^2=0 \}$.

\no {\rm (iv)} Beam equation: we have $\partial_x g_2(t,dx)=G'({x \over \sqrt{t}})dx$.

\end{Example}

\no In order to introduce the probabilistic representation of the solution of \eqref{IVP2}-\eqref{BC2}, we consider the branching mechanism defined in Subsection \ref{sect:branching}, where we modify the definition of the independent iid random variables $(J_k)_{k\in\mathbf{K}}$ and we introduce the types of particles $(\theta^t_k)_{k\in\mathbf{K}}$ as follows:
\begin{itemize}
\item $\P[J_k=\ell_j]=q_j$ for all $j\ge 0$, and we denote $\bar J_{k,-1}:=0$, $\bar J_{k,h}:=\bar J_{k,h-1}+J_{k,h}$, $h=0,\ldots,H$;
\item an arbitrary particle $k\in\overline{\Kc}_t\setminus\Kc_t$ branches at time $T_k$ into $\bar J_{k,H}$ particles; the $h-$th block of  descendantparticles are labelled
 \be\label{k-j}
 (k,j),~~j=\bar J_{k,h-1}+1,\ldots, \bar J_{k,h}
 &h=0,\ldots,H;&
 \ee
\item we assign to the $h-$th block of new particles \eqref{k-j} the types:
 \begin{equation}
 \theta^t_{(0)}:=0
 ~\mbox{and}~
 \theta^t_{(k,j)} := h,
 ~\mbox{for}~~k\in\overline{\Kc}_t\setminus\Kc_t,~j=\bar J_{k,h-1}+1,\ldots, \bar J_{k,h},~h=0,\ldots,H.
 \end{equation}
\end{itemize}

\no Finally, we introduce a branching process which differs slightly from \eqref{Xk}-\eqref{Zk}. Let
 \be\label{hatXk}
 \hat X^0_0=x,
 &\mbox{and}&
 \hat X^k_{T^t_k}
 :=
 \hat X^{k_-}_{T^t_{k_-}} + \hat Z^k_{T^t_k},
 ~~\mbox{for all}~~
 k\in\overline{\Kc}_t,
 \ee
where the distribution of $\hat Z^k_{T^t_k}$ depends on the type of the particle:
 \be\label{hatZk}
 \P\big[ \hat Z^k_{T^t_k}\in dz \big| I^t_k,\theta^t_k, \Delta T^t_k \big]
 &=&
 \1_{\{\theta^t_k=0\}}\mu_{I^t_k}(\Delta T^t_k, dz)
 +\1_{\{\theta^t_k\neq 0\}}\mu^1_{I^t_k}(\Delta T^t_k, dz),
 \ee
with $\Delta T^t_k:=T^t_k-T^t_{k_-}$. The main goal of this section is to provide a representation of the solution of the Cauchy problem \eqref{IVP2}-\eqref{BC2} by means of the random variable
 \begin{equation}
 \hat\xi_{t,x}
 :=
 \prod_{k\in\Kc_t}  \frac{\Wc_k}
                                     {\bar\rho\big(\Delta T^t_k\big)}
            \big[f_{I^t_k}\big(\hat X^{k}_t\big)-\1_{\{\theta^t_k\neq 0\}}f_{I^t_k}\big(\hat X^{k}_{T_{k-}}\big)\big]
 \;       \prod_{k\in\overline{\Kc}_t\setminus \Kc_t} \frac{\Wc_k}
                                                                                      {\rho\big(\Delta T^t_k\big)}
                                                                    c_{J_k,0}\big(t-T^t_k,\hat X^k_{T^t_k}\big),
 \label{xi-semilinear2}
 \end{equation}
where the random weights $\Wc_k$ are given by
 \b*
 \Wc_k
 :=
 \1_{\{\theta^t_k=0\}}\gamma_{I^t_k}\big(\Delta T^t_k,\hat Z^k_{T^t_k}\big)
 + \1_{\{\theta^t_k\neq 0\}}c_{J_k,\theta^t_k}\big(t-T^t_{k-},\hat X^{k-}_{T^t_{k-}}\big)
                                         \cdot
                                         \gamma^1_{I^t_k}\big(\Delta T^t_k,\hat Z^k_{T^t_k}\big),
 &k\in\overline{\Kc}_t.&
 \e*

\no Recall the power series $H_p$ introduced in \eqref{cond:Fp}, and define
 \be
 \!\!\!\!\!\!\!\!\hat r_p
 &\!\!\!\!\!\!\!\!:=&
 \!\!\!\!\!\!\!\!\!\!\!\!
 \sup_{\tiny\begin{array}{c}
           0\!\le\! t\!\le\! T\\1\!\le\! n\!\le\! N
           \end{array}}
           \!\!\!\!
 \frac{\big\{\|f_n\|_\infty^p \int |\gamma_n(t,z)|^{p-1}|g_n|(t,dz)\big\}
          \vee
          \big\{\|\nabla f_n\|_\infty^p \int |z|^p|\gamma^1_n(t,z)|^{p-1}|\partial_xg_n|(t,dz)\big\}
          }
        {\bar\rho(t)^{p-1}}
        ~~
 \label{Ap}
 \\
 \hat\alpha_p
 &:=&
 \sup_{0\le t\le T}
 \frac{\big\{\int |\gamma_N(t,z)|^{p-1}|g_N(t,dz)|\big\}
          \vee
          \big\{\int |\gamma^1_N(t,z)|^{p-1}|\partial_xg_N(t,dz)|\big\}
          }
        {\rho(t)^{p-1}}
 .
 \label{Bp}
 \ee
The following result provides a probabilistic representation of the solution of the semilinear Cauchy \eqref{IVP2} under a condition involving the above $\hat r_p$ and $\hat\alpha_p$, that will be guaranteed in Proposition \ref{prop:condIVP2} below to hold under some sufficient conditions.

 \begin{Theorem}\label{thm:IVP2}
Let Assumptions \ref{hyp:hatg}, \ref{hyp:g}, and \ref{hyp:Dgn} hold true, and assume that $f$ is bounded and Lipschitz. Let $p>1$, and $\rho$ be a positive density function with support on $(0,\infty)$ such that the constants $\hat r_p$ and $\hat \alpha_p$ defined in \eqref{Ap}-\eqref{Bp}  satisfy:
 \be\label{condIVP2}
 \hat r_p<R_p, ~\hat\alpha_p<\infty
 &\mbox{and}&
 \int_{\hat r_p}^{\hat s_p} \frac{ds}{H_p(s)} = \hat\alpha_p T
 ~\mbox{for some}~\hat s_p > \hat r_p~\mbox{and}~T>0.
 \ee
Then, $\hat\xi_{t,x}\in\L^p$ for $t\in[0,T]$, and  $u(t,x):=\E\big[\hat\xi_{t,x}\big]$ is the unique bounded continuous solution of the semilinear Cauchy problem \eqref{IVP2}.
\end{Theorem}

\proof Similar to the proof of Theorem \ref{thm:IVP1}, we directly estimate that
 \b*
 \big|\hat\xi_{t,x}\big|^p
 &\le&
 \prod_{k\in\Kc_t} \frac{1}{\bar\rho(\Delta T_k^t)^p}
                             \Big(\1_{\theta^t_k=0}\big\|f_{I_k^t}\big\|_\infty^p
                                  \big|\gamma_{I^t_k}(\Delta T_k^t,\hat Z^k_{T_k^t})\big|^p
                                    + \1_{\theta^t_k\neq 0}\big\|\nabla f_{I_k^t}\big\|_\infty^p
                                                                  \big|Z^k_{T_k^t}\big|^p
                                                                  |\Wc_k|^p
                             \Big)
 \\
 &&\hspace{10mm}\times
 \prod_{k\in\overline{\Kc}_t\setminus\Kc_t} \Big|\frac{\Wc_k}{\rho(\Delta T^t_k)}\Big|^p
                                                                   \|c_{J_k,0}\|_\infty^p.
 \e*
By the independence of the $\Delta T_k^t$'s and the $Z^k_{T_k^t}$'s, this provides
 \b*
 \E\big|\hat\xi_{t,x}\big|^p
 &\le&
 \E\prod_{k\in\Kc_t} \frac{\hat r_p}{\bar\rho(\Delta T_k^t)}
    \prod_{k\in\overline{\Kc}_t\setminus\Kc_t} \frac{\hat\alpha_p\,\|c_{J_k,0}\|_\infty^p}{\rho(\Delta T_k^t)}.
 \e*
The required result follows by the same line of argument as in the proof of Theorem \ref{thm:IVP1}.
\ep


\no In the rest of this section, we provide sufficient conditions which guarantee that the constants $\hat r_p$ and $\hat\alpha_p$, involved in Condition \eqref{condIVP2}, are finite. In preparation for this, we need some estimates on the Green functions $g_n$, $n=1,\ldots,N$. Recall that $\{\alpha\in\N^d_M:|\alpha|=M~\mbox{and}~b_\alpha\neq 0\}\neq\emptyset$, and define {\it the principal symbol}
 $$
 b_M(\xi)
 :=
 i^{ M}\!\!\sum_{|\alpha|=M} b_\alpha \xi^\alpha
 =
 e^{i\pi\eta_M(\xi)}|b_M(\xi)|,
 ~\mbox{where}~
 \eta_M(\xi):=\frac{1+M-\mbox{sg}\{b_M(\xi)\}}{2};
 ~\xi\in\R^d.
 $$

\begin{Lemma}\label{lem:lambdaeps}
For all $\xi\in\R^d$, the matrix $B(\eps^{-1}\xi)$ has $N$ simple eigenvalues $\lambda_n(\eps^{-1}\xi)$, $n=1,\ldots,N$, for sufficiently small $\eps>0$, with asymptotics
 \b*
 \lim_{\eps\searrow 0}\eps^{\frac{1}{\sigma}}\lambda_n\big(\eps^{-1}\xi\big)
 =
 \lambda^0_n(\xi)
 :=
 \big|b_M(\xi)\big|^{\frac{1}{N}}
 e^{\frac{i\pi}{N}(\eta_M(\xi)+2n)},
 &\mbox{where}~~
 \sigma := \frac{N}{M}.
 &
 \e*
\end{Lemma}

\proof Recall from Remark \ref{rem:CharacteristicPolynomial} that the spectrum of the matrix $B\big(\eps^{-1}\xi\big)$ consists of the solutions of the characteristic polynomial $ \sum_{n=1}^N a_n \lambda^n = b(\eps^{-1}\xi)$.

\no As $\{\alpha\in\N^d_M:|\alpha|=M~\mbox{and}~b_\alpha\neq 0\}\neq\emptyset$, we see that $\eps^M b(\eps^{-1}\xi\big)\longrightarrow b_M(\xi)$. Then, denoting by $\lambda_\eps$ an arbitrary solution of the last characteristic polynomial, we deduce that $(\lambda_\eps)_{\eps>0}$ has no finite accumulation point. Together with the normalization $a_N=1$, this in turn implies that $\sum_{n=1}^N a_n \lambda_\eps^n\sim \lambda_\eps^N$, and therefore
 \b*
 \lim_{\eps\searrow 0} \eps^M\lambda_\eps^N
 &=&
 b_M(\xi)
 \;=\;
 \big|b_M(\xi)\big| e^{i\pi \eta_M(\xi)}
 \;:=\;
 \lambda^0_N(\xi)^N,
 \e*
which can be written equivalently as $\lim_{\eps\searrow 0} \big(\eps^{\frac{M}{N}} \frac{\lambda_\eps}{\lambda^0_N(\xi)}\big)^N=1$. Hence, the limiting spectrum of the matrix $B\big(\eps^{-1}\xi\big)$ consists of $N$ simple eigenvalues expressed in terms of the unit roots:
 \b*
 \lim_{\eps\searrow 0}\eps^{\frac{M}{N}}\lambda_n\big(\eps^{-1}\xi\big)
 &=&
 \lambda^0_N(\xi) e^{\frac{2ni\pi}{N}}
 \;=\;
 \lambda^0_n(\xi),
 ~~n=1,\ldots,N.
 \e*
\ep

\no The last result shows that the asymptotics of $B(\eps^{-1}\xi)$ are in the context of Remark \ref{rem:CharacteristicPolynomial}. Then, it follows that the space Fourier transforms of the Green functions are given by \eqref{hatgn:rem}. Define the corresponding limits:
 \begin{eqnarray}\label{hatgn0}
 \hat g^0_n
 :=
 \Lambda^0_n\sum_{j=1}^N \frac{e^{\lambda^0_j}}
 {\prod_{\ell\neq j}(\lambda^0_\ell-\lambda^0_j)},
 &n=1,\ldots,N,&
 \end{eqnarray}
with:
 $$
 \Lambda^0_N(\xi):=1,
 ~~
 \Lambda^0_n(\xi)
 :=
 (-1)^{n-1}\!\!\!\!\sum_{_{\!\!\!\tiny\begin{array}{c}
                    1\!\!\le\!\!\ell_1\!\!\le\!\!\ldots\!\!\le\!\!\ell_{N-n}\!\!\le\!\! N
                                                                                      \\
                                                                                      \ell_1,\ldots,\ell_{N-n}\!\!\neq\!\! n
                    \end{array}}}
                    \!\!\!\!(\lambda^0_{\ell_1}\cdots\lambda^0_{\ell_{N-n}})(\xi)
 ~~\mbox{for}~n<N.
 $$
The following additional condition is needed in order to characterize the short time asymptotics of the Green functions.

 \begin{Assumption}\label{ass:UI}
 For all $\varphi\in\Sc$, the family $\{t^{1-n}\hat g_n(t,t^{-\sigma}\cdot)\varphi, t\in(0,\eps]\}$ is uniformly integrable in $\L^1(\R^d)$, for some $\eps>0$.
 \end{Assumption}

\no We denote by $X$ the canonical map on $\R^d$, i.e., $X(x)=x$  for all $x\in\R^d$, and we introduce the scaled Green functions
 \b*
 g^\sigma_n(t,\cdot) := g_n(t,\cdot)\circ(t^\sigma X)^{-1},
 &t\in[0,T],&
 n=1,\ldots,N.
 \e*
In particular, in the case where $g_n$ can be represented by a function, $g^\sigma_n(t,x)=t^{\sigma d}g_n(t,t^\sigma x)$.

\begin{Lemma}\label{lem:gnwn}
The space Fourier Green functions $\hat g$ satisfy the short time asymptotics
 \b*
 t^{1-n}\hat g_n(t, t^{-\sigma}\xi)
 \;\longrightarrow\;
 \hat g_N^0(\xi)
 &\mbox{as $t\searrow 0$, for all}&
 \xi\in\R^d.
 \e*
If in addition Assumption \ref{ass:UI} holds true, then the last convergence holds in $\Sc'$, and the short time asymptotics of the scaled Green functions are given by
  \b*
  t^{1-n} g^\sigma_n(t,\cdot)
  \longrightarrow
  g_n^0 := \mathfrak{F}^{-1}\hat g^0_n
  &\mbox{as $t\searrow 0$, in}&
  \Sc',
  ~~n=1,\ldots,N.
  \e*
\end{Lemma}

 \proof
First, the pointwise convergence of $\hat G_n(t,\cdot):=t^{1-n}\hat g_n(t,t^{-\sigma}\cdot)$ towards $\hat g_n^0$, as $t\searrow 0$, follows from direct application of Lemma \ref{lem:lambdaeps} together with the observations reported in Remark \ref{rem:CharacteristicPolynomial}. The uniform integrability condition of Assumption \ref{ass:UI} guarantees that $\langle\hat G_n(t,\cdot),\varphi\rangle \longrightarrow\langle\hat g_n^0(t,\cdot),\varphi\rangle$ for all $\varphi\in\Sc$, i.e., $\hat G_n(t,\cdot)\longrightarrow\hat g^0_n$ as $t\searrow 0$ in $\Sc'$. This in turn implies the convergence of the corresponding Fourier inverse $\mathfrak{F}^{-1}\hat G_n(t,\cdot)$ towards $g_n^0 := \mathfrak{F}^{-1}\hat g^0_n$ as $t\searrow 0$ in $\Sc'$. It remains to relate the distribution $\mathfrak{F}^{-1}\hat G_n(t,\cdot)$ to the Green function $g_n$. To see this, we use the properties of the Fourier transform in $\Sc'$ as defined by means of arbitrary test functions $\varphi\in\Sc$ as follows:
 \b*
 \langle \mathfrak{F}^{-1}\hat G_n,\varphi \rangle
 =
 \langle \hat G_n,\mathfrak{F}^{-1}\varphi \rangle
 =
 t^{\sigma d}
 \langle \hat g_n(t,\cdot),(\mathfrak{F}^{-1}\varphi)(t^\sigma\cdot) \rangle
 =
 t^{\sigma d}\langle g_n(t,\cdot),\mathfrak{F}\big((\mathfrak{F}^{-1}\varphi)(t^\sigma\cdot)\big) \rangle.
 \e*
We finally observe by direct calculation that $\mathfrak{F}\big((\mathfrak{F}^{-1}\varphi)(\lambda\cdot)\big)=\lambda^{-d}\varphi(\lambda^{-1}\cdot)$ for all $\varphi\in\Sc$ and all constant $\lambda\in\R$. Then
 \b*
 \langle \mathfrak{F}^{-1}\hat G_n,\varphi \rangle
 =
 \langle g_n(t,\cdot),\varphi(t^{-\sigma}\cdot) \rangle
 =
 \langle g^\sigma_n(t,\cdot),\varphi \rangle.
 \e*
By the arbitrariness of $\varphi\in\Sc$, this provides that $\mathfrak{F}^{-1}\hat G_n(t,\cdot)=g^\sigma_n(t,\cdot)$, thus completing the proof.
 \ep


\no We observe that, as a consequence of the convergence of the scaled Green functions in $\Sc'$, we deduce that
 \be\label{convDalphagn}
 t^{1-n+\sigma |\alpha|} \partial_x^\alpha g^\sigma_n(t,\cdot)
  \longrightarrow
  \partial_x^\alpha g_n^0
  &\mbox{as $t\searrow 0$, in}&
  \Sc',
  ~~n=1,\ldots,N,
 \ee
where $\partial_x^\alpha:=\frac{\partial^{|\alpha|}}{\partial x_1^{\alpha_1}\ldots\partial x_d^{\alpha_d}}$ for all $\alpha\in\N^d$. Our final result requires the following conditions.

\begin{Assumption}\label{hyp:rsfinis}
{\rm (i)} The functions $t\longmapsto \int |z|^p|\gamma^1_n(t,z)|^p\mu_n^1(t,dz)$, $n=1,\ldots,N,$ and $t\longmapsto \int |\gamma^1_N(t,z)|^p\mu_N^1(t,dz)$ are continuous on $(0,T]$;
\\
{\rm (ii)} $\gamma_n^1(t,t^\sigma z)={\rm O}(1)$ near the origin $t=0$;
\\
{\rm (iii)} The families
 \b*
 &\big\{z\longmapsto |z|^p |\gamma_n^1(t,t^\sigma z)|^{p-1}|\partial_xg^\sigma_n|(t,dz)\big\}_{t\le\eps},~1\le n\le N,&
 \\
 &\mbox{and}~
 \big\{z\longmapsto |\gamma_N^1(t,t^\sigma z)|^{p-1}|\partial_xg^\sigma _N|(t,dz)\big\}_{t\le\eps}&
 \e*
are uniformly integrable for some $\eps>0$.
\end{Assumption}

\begin{Proposition}\label{prop:condIVP2}
Let Assumptions \ref{hyp:hatg}, \ref{hyp:g}, \ref{hyp:Dgn} (i), \ref{ass:UI} and \ref{hyp:rsfinis} hold true. Assume further that the density $\rho\in C^0(\R_+)$, strictly positive on $(0,\infty)$, and $\lim_{t\searrow 0} t^{N-1+\sigma p}\rho(t)^{1-p}<\infty$. Then, for any bounded Lipschitz function $f$, we have $\hat r_p+\hat\alpha_p<\infty$.
\end{Proposition}

\proof By Assumption \ref{hyp:Dgn} (i), together with the fact that $\bar\rho(T)\le\rho\le 1$ on $[0,T]$, we only need to justify that
 \b*
 \sup_{\tiny\begin{array}{c}
           0\!\le\! t\!\le\! T\\1\!\le\! n\!\le\! N
           \end{array}}
           \!\!\!\!
 \int |z|^p|\gamma^1_n(t,z)|^{p-1}|\partial_xg_n|(t,dz)\big)
 +
 \sup_{0\le t\le T}
 \rho(t)^{1-p}\int |\gamma^1_N(t,z)|^{p-1}|\partial_xg_N|(t,dz)\big)
 <\infty.
 \e*
\no Assumption \ref{hyp:rsfinis} (i) ensures that the functions inside the last suprema are continuous on $(0,\infty)$. Then, in order to prove that $\hat r_p$ and $\hat\alpha_p$ are finite, it suffices to verify that the functions inside the last suprema are bounded near $t=0$.

\no To see this, we compute by a direct change of variables that:
 \b*
 \int |z|^p|\gamma^1_n(t,z)|^{p-1}|\partial_xg_n|(t,dz)
 &=&
 \int |t^\sigma z|^p |\gamma^1_n(t,t^\sigma z)|^{p-1}|\partial_xg^\sigma_n|(t,dz)
 \\
 &=&
 t^{n-1+\sigma p}
 \int |z|^p t^{1-n}|\gamma^1_n(t,t^\sigma z)|^{p-1}|\partial_xg^\sigma_n|(t,dz)
 \\
 &=&
 \mbox{O}\big(t^{n-1+\sigma p}\big),~~\mbox{near the origin},
 \e*
by Assumption \ref{hyp:rsfinis} (ii)-(iii), together with the consequence \eqref{convDalphagn} of Lemma \ref{lem:gnwn}. This implies that the limit is zero as $n\ge 1$ and $p>1$. Similarly,
 \b*
 \rho(t)^{1-p}\int |\gamma^1_N(t,z)|^{p-1}|\partial_xg_N|(t,dz)
 &=&
 t^{N-1+\sigma p}\rho(t)^{1-p}
 \int t^{1-N}|\gamma^1_N(t,t^\sigma z)|^{p-1}|\partial_xg^\sigma_N|(t,dz)
 \\
 &=&
 \mbox{O}\big(t^{N-1+\sigma p}\rho(t)^{1-p}\big),~~\mbox{near the origin},
 \e*
again by Assumption \ref{hyp:rsfinis} (ii)-(iii).
\ep

\section{Numerical examples}
\label{sect:numerics}

\subsection{Wave semi-linear PDE} \label{Wave semi-linear PDE}

We consider the nonlinear Klein-Gordon wave equation in $\RR^d$ for $1 \leq d\leq 3$:
 \bea\label{nlKG}
 (\partial^2_{tt} - \Delta) u+u^3+u^2
 =
 0,
 \quad u(0,x)=f_1(x),
 \quad  \partial_t u(0,x)=f_2(x).
 \eea
To the best of our knowledge, the current literature only considers approximation of the solution by deterministic numerical schemes, see e.g. \cite{deh}. Due to the curse of dimensionality, mainly $d=1$ and $d=2$ have been considered. To illustrate the efficiently of our algorithm, we solve this equation in $d=1,2$ and $3$. In our numerical experiments, we take the initial conditions
 \beaa
 f_1(x):=-{12 \over 9+2(\sum_{i=1} ^d x_i)^2},
 &\mbox{and}&
 f_2(x):= -\frac{48\sqrt{d+1}  (\sum_{i=1} ^d x_i)}{\left(2 (\sum_{i=1} ^d x_i)^2+9\right)^2},
 \eeaa
 for which the explicit solution is $u(t,x)=-{12 \over 9+2( \sqrt{d+1}t - \sum_{i=1} ^d x_i)^2}$.

\no We choose $\rho(t)=\beta e^{-\beta t}$ and $\bar{\rho}(t)=e^{-\beta t}$. For convenience, we  set  $u(t,x):=U(t,x)+f_1(x)$, and we compute that $U$ satisfies the non-homogeneous nonlinear wave PDE:
 \beaa
 &&
 (\partial^2_{tt} - \Delta) U+U^3+(3f_1+1)U^2+(3f_1^2 +2f_1)U+(f_1^3+f_1^2- \Delta f_1)=0
 \\
 &&U(0,x)=0, \quad  \partial_t U(0,x)=f_2(x).
 \eeaa
Note that as $U(0,x)=0$, we do not need to simulate our branching particles according to the distribution $g_1$ but only $g_2$. This was our motivation for introducing the function $U$.

\no We directly compute that $r_p< ({1+d \over 6})^p t^p$ and $\alpha_p<t^p \beta^{1-p}$. Furthermore, $|(3f_1+1)|_\infty=1$,
$|(3f_1^2 +2f_1)|_\infty=8/3$ and $|(f_1^3+f_1^2- \Delta f_1)|_\infty<1$. This implies that the assumptions in Theorem \ref{Thmsemilinear} are satisfied.

\no  The Monte-Carlo approximation of $U(t,X_0)$ (and therefore $u(t,X_0)=U(t,X_0)+f_1(X_0)$) can be described by the following meta-algorithm:

\subsubsection*{Meta-algorithm:}
\begin{enumerate}
\item Start at $t_0=0$ at the position $X_0$ and initialize a weight ${\cal W}:=1$.
\item Simulate an exponential r.v. $\tau$ with (arbitrarily) constant intensity $\beta$ and simulate the particle at the new position $X_\tau=X_0+Z_\tau$ at $\tau$. More precisely, we draw uniform variables $(U_i)_{i=1,2}$ on $[0,1]$ and set
 $$
 \begin{array}{l}
 Z^1_{\tau}= \tau (2U_1-1), ~d=1
 \\
 Z_{\tau}^1=\sqrt{1-U_1^2}\cos(2\pi U_2)\tau,
 ~Z_{\tau}^2=\sqrt{1-U_1^2}\sin(2\pi U_2)\tau,~d=2
 \\
 Z_{\tau}^1=\cos(2\pi U_1)\cos(2\pi U_2)\tau,
 ~Z_{\tau}^2=\sin(2\pi U_1)\cos(2\pi U_2)\tau,
 ~Z_{\tau}^3=\sin(2\pi U_2)\tau, ~d=3.
 \end{array}
 $$
At $t_1:=\tau<t$, the particle dies and we create $0$, $1$, $2$ or $3$ descendants with probability $p:=1/4$. We then  multiply the weight ${\cal W}$ by the mass $\tau$ and according to the number of descendants, we update again the weight ${\cal W}$ by
 \beaa
 {\cal W}
 &:=&
 \left\{\begin{array}{l}
         {\cal W} \times (-p^{-1})(f_1^3+f_1^2- \Delta f_1)(X_\tau)\beta^{-1}e^{\beta(t_1-t_{0})},~
         \mbox{if}~0 \; \mathrm{descendant,}
         \\
         {\cal W} \times (-p^{-1})(3f_1^2 +2f_1)(X_\tau)\beta^{-1}e^{\beta(t_1-t_{0})},~
         \mbox{if}~1\; \mathrm{descendant}
         \\
         {\cal W} \times (-p^{-1})(3f_1+1)(X_\tau)\beta^{-1}e^{\beta(t_1-t_{0})},~
         \mbox{if}~2\; \mathrm{descendants.}
         \\
         {\cal W} \times (-p^{-1})\beta^{-1}e^{\beta(t_1-t_{0})},~
         \mbox{if}~3\; \mathrm{descendants,}
         \end{array}
 \right.
 \eeaa
where $0$ descendant means that the particle dies.
\item For each descendant, we apply independently Steps 2 and 3 until the default time -- say $\tau_n$ -- is greater than the maturity $t$.  In this case, we multiply ${\cal W}$ by
 \beaa
 {\cal W}
 &:=&
 {\cal W} \times  e^{\beta(t-\tau_{n-1})}
 \eeaa
\item Finally, for all particles alive at time $t$ (with locations $(X_t^k)_{k\in\Kc_t}$), compute
 \beaa
 {\cal W} \prod_{k\in\Kc_t} f_2(X^k_t),
 \eeaa
and average the result using $M$ Monte-Carlo paths.
\end{enumerate}

\begin{figure}[h]
\begin{center}
\includegraphics[width=8cm,height=6cm]{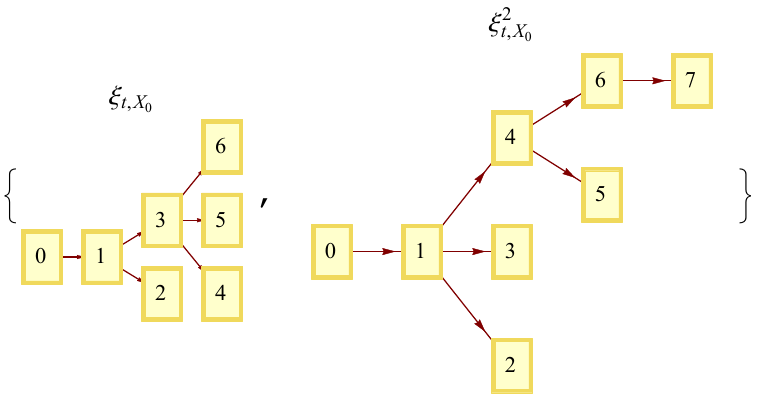}
\end{center}
\caption{Galton-Watson trees associated to $ \xi^{(1)}_{t,X_0}$ and $ \xi^{(2)}_{t,X_0}$.}
\label{FigGW}
\end{figure}

\no In Figure  \ref{FigGW}, we  have plotted two examples of Galton-Watson trees corresponding to the  functionals: 
 \beaa
 \xi^{(1)}_{t,X_0}
 &\!\!\!\!:=&\!\!\!\!
 {1 \over  \beta^2 p^2}(3f_1(X_{\tau_1}+1)e^{\beta \tau_1}
 e^{\beta (\tau_3-\tau_1)}e^{3\beta (t-\tau_3)}e^{\beta (t-\tau_1)}f_2(X_t^6)f_2(X_t^5)f_2(X_t^4)f_2(X_t^2),
 \\
 \xi^{(2)}_{t,X_0}
 &\!\!\!\!:=&\!\!\!\!
 -\frac{1}{\beta^3 p^3}(3f_1(X_{\tau_4}+1)
 (3f_1^2+2f_1)(X_t^6)e^{\beta \tau_1}
{e^{\beta (\tau_4-\tau_1)}}
 {e^{\beta (\tau_6-\tau_4)}}
 f_2(X_t^7)f_2(X_t^5)f_2(X_t^3)f_2(X_t^2)\\
 &&e^{\beta(t-\tau_6)}e^{\beta(t-\tau_4)}e^{2\beta(t-\tau_1)}.
 \eeaa
In the present very simple examples, in order to alleviate the figure, we have used simpler notations to label the branching particles than those in Section \ref{sect:branching}.

\begin{figure}[h]
\begin{center}
\includegraphics[width=6cm,height=6cm]{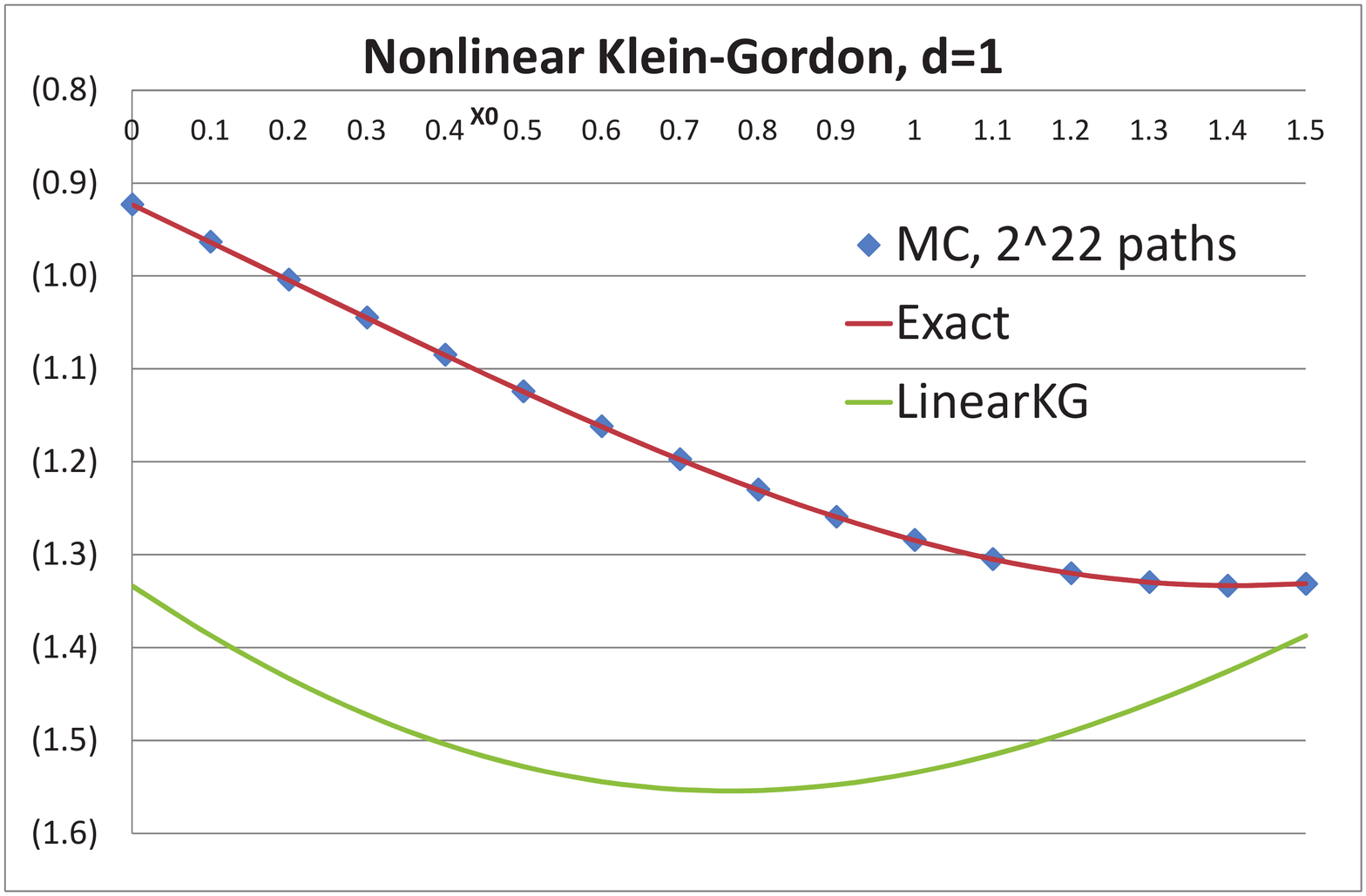}
~~~~~~~~
\includegraphics[width=6cm,height=6cm]{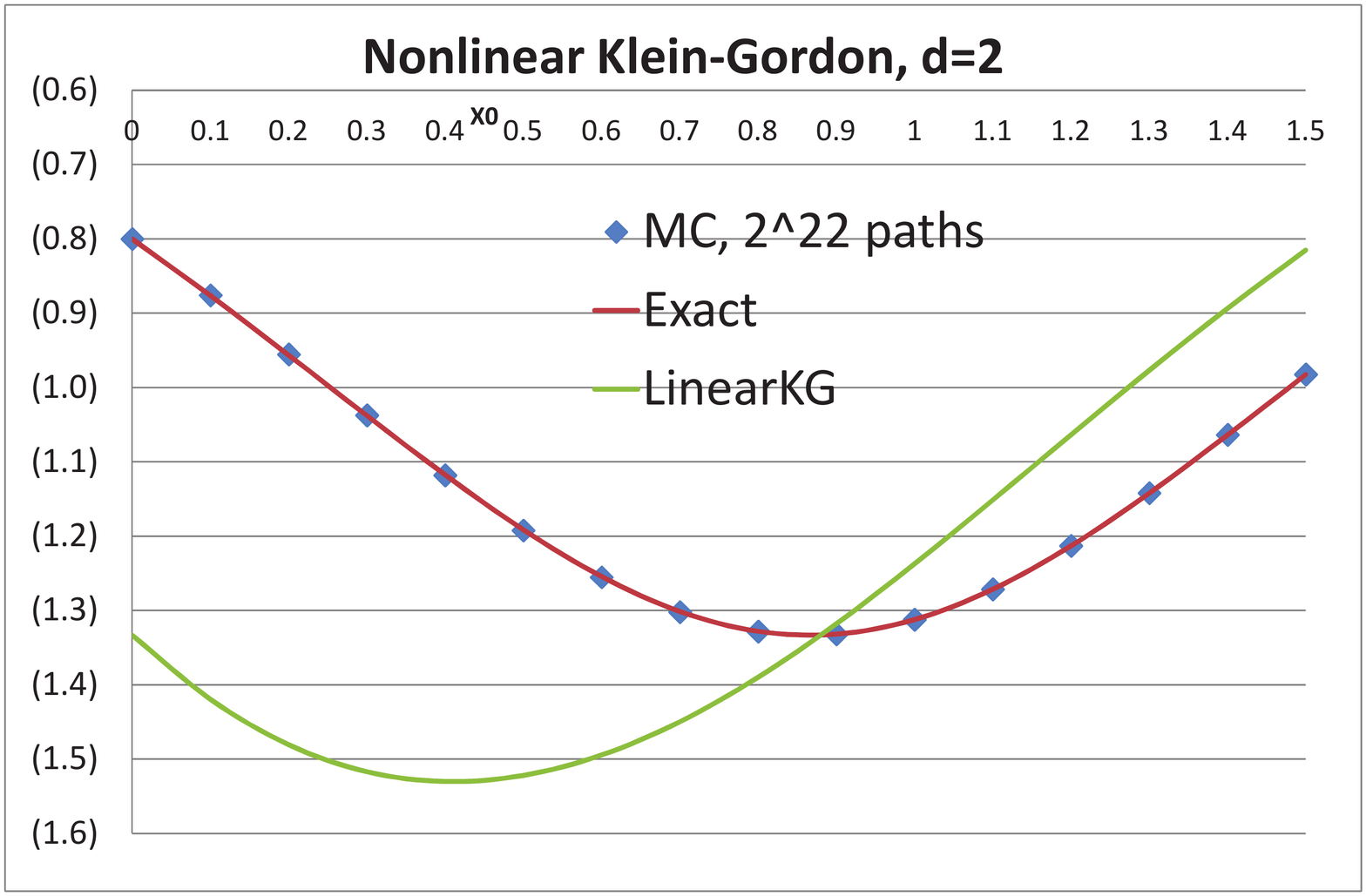}
\includegraphics[width=6cm,height=6cm]{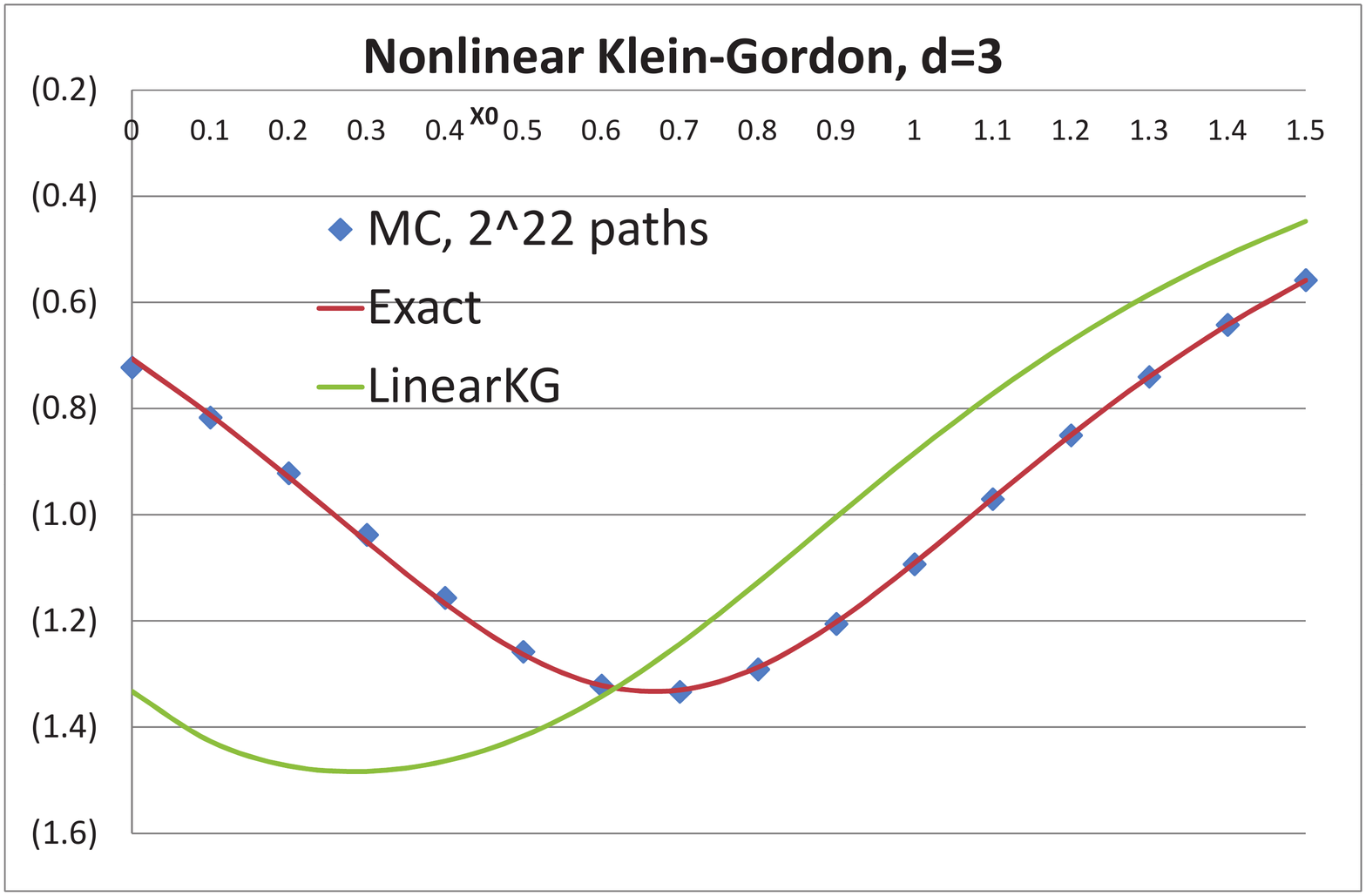}
\end{center}
\caption{Numerical solutions of the nonlinear Klein-Gordon PDE (\ref{nlKG}) for $d=1,2,3$ as a function of $x_0$.}
\label{Fig2}
\end{figure}

\no Note that by construction, the result is independent of $\beta$ when $M$ is large. There is an optimal choice of $\beta$ that minimizes the variance of our Monte-Carlo estimator. We have chosen $\beta=1$ in our numerical experiments and $M=2^{22}$ for which the standard deviation of our estimator is less than $0.01$. Below, we have plotted our numerical result for $u(t,X_0)$ again our analytical solution for $d=1, 2, 3$ as a function of $x_0 \in [0,1.5]$ (all the coordinates in $\RR^d$ are equal to $x_0$) and $t=1$, see Figure \ref{Fig2}.  We obtain a perfect match. In order to see that our numerical solutions captures perfectly the additional nonlinearity $u^3+u^2$, we have also shown  for completeness the (analytical) solution of the linear wave equation (denoted ``LinearKG''):
 \bea\label{lKG}
 (\partial_{tt} - \Delta) u=0, \quad u(0,x)=f_1(x), \quad  \partial_t u(0,x)=f_2(x).
 \eea

\subsection{Yang-Mills PDE: a toy model}
\label{sect:YangMills}

\no Here we consider the semilinear  wave equation in $\RR$:
 \bea\label{scalarYM}
 \partial_{tt}u -\Delta u + u^3 +u \partial_{x} u=0, \quad u(0,x)=f_1(x), \quad  \partial_t u(0,x)=f_2(x).
 \eea
We take the initial conditions
 \beaa
 f_1(x)=-\frac{1}{1 -x}, \quad f_2(x)=\frac{1}{(1 - x)^2},
 \eeaa
for which the explicit solution is $ u(t,x)=-\frac{1}{1+  t - x} $.   Assumption (\ref{hyp:g}) is satisfied for $g_2$ only for $d=1$, this is why we restrict to this case. For $d>1$, $\nabla g_2$ involves derivative of delta function supported on the lightcone. Notice that the singularity of $f_1$ and $f_2$ are not seen by  our numerical algorithm because of our choice of the initial position $x_0$ and the finite speed property of the wave equation: The singularity $x=1$ is not attainable.

\no This example can be interpreted as a scalar version of the Yang-Mills hyperbolic system.  Indeed, the Lagrangian associated to the Yang-Mills theory on $\RR^{1+d}$ is
 \beaa
 \int  F_{\mu \nu}^a  F^{\mu \nu a} d^4 x, \quad  F^{\mu \nu a}:=\eta^{\mu \alpha} \eta^{\nu \beta} F_{\alpha \beta}^ a
 \eeaa with $\eta_{00}=-1$, $\eta_{ii}=1$, $i=1,\cdots d$ and $0$ otherwise. Here we  use the convention of implicit summation for repeated indices $(a,\mu,\nu)$. The curvature is $F_{\mu \nu}^a:=\partial_\mu A_\nu^a -\partial_\nu A_\mu^a -f_{bc}^a A_\mu^b A_\nu^c$. The Euler-Lagrange equations (written in the gauge frame $\partial_\mu A^{\mu,a}=0$) give the  system of hyperbolic PDE
\beaa \partial_\mu F^{\mu \nu a}+f_{bc}^a A_\mu^b F^{\mu \nu c}=0  \eeaa
Note that the initial boundary conditions at $t=0$, $A_{\nu}^a(t=0,x)$ and $\partial_t A_{\nu}^a(t=0,x)$ need to satisfy a constraint condition in order to preserve the gauge condition for all $t$.  Using the expression for the curvature $F_{\mu \nu}^a$, we get that $A^{\mu a}:=\eta^{\mu \alpha } A_{\alpha}^a $ is solution of  a system of hyperbolic PDE of the form:
\beaa \square A ^{\nu a} -f^a_{bc} A^{\mu b} \partial_\mu A^{\nu c}+
f^{a}_{bc} A_\mu^b \left(\partial^\mu A^{\nu c} -\partial^\nu A^{\mu c} -f_{de}^c A^{\mu d} A^{\nu e} \right)=0  \eeaa
\no with $\square:=\partial_{tt}-\Delta$. This can schematically written as
\beaa \square A - f A \partial A - f^2 A^3 =0 \eeaa  hence our PDE (\ref{scalarYM}).

\no Note that our branching Monte-Carlo algorithm can be easily adapted to solve 	a system of semilinear PDEs. All we need to do is to index particles with a type $(\mu a)$ corresponding to a coordinate of the solution, see an example with the complex Gross-Pitaevskii equation in Section \ref{Gross--Pitaevskii}.

\no As in  the previous section, we set $u(t,x)=f_1(x)+U(t,x)$ for which $U$ satisfies the non-homogeneous hyperbolic PDE
 \beaa
 \partial_{tt}U -\Delta U + U^3 +3 f_1 U^2+ U \sum_{i=1}^d \partial_{x_i} U+(3 f_1^2(x)  +   \sum_{i}^d  \partial_{x_i} f_1)U  +
 f_1 \sum_{i}^d  \partial_{x_i} U \\+(f_1^3+\sum_{i}^d f_1 \partial_{x_i} f_1-\Delta f_1)=0, \quad U(0,x)=0, \quad  \partial_t u(0,x)=f_2(x)
 \eeaa
\no  The Monte-Carlo approximation of $U(t,X_0)$ (and therefore $u(t,X_0)=U(t,X_0)+f_1(X_0)$) can be described by the following meta-algorithm:

\begin{enumerate}
\item Start with a type $(0)$ particle at $t_0=0$ at position $X_0$, and initialize a weight ${\cal W}:=1$;
\item Simulate an exponential r.v. $\tau$ with (arbitrarily) constant intensity $\beta$, and simulate the particle at the new position $X_\tau=X_0+Z_\tau$ at $\tau$. More precisely, we draw an uniform variable $(U)$ on $[0,1]$ and set $Z_{\tau}= \tau (2U-1)$;
\item At $t_1:=\tau<t$, the particle dies and we create $0$, $1$, $2$ or $3$ descendants with probability $p:=1/6$; in case of $1$ descendant, the type assigned is $(0)$ or $(1)$; similarly, in case of $2$ descendants, their types can be both $(0)$, both $(1)$, or one of each type; in case of $3$ descendants, they are all of type $(0)$; we then  multiply the weight ${\cal W}$ by the mass $\tau$ (see Example \ref{eg:mu}-(i)) if the type of the particle is $(0)$ and according to the number of descendants and type, we update also the weight ${\cal W}$ by
 \beaa
 {\cal W}
 &:=&
 \left\{\begin{array}{l}
 {\cal W} -p^{-1}(f_1^3+ f_1 \partial_{x} f_1-\Delta f_1)(X_\tau)\beta^{-1}e^{\beta(t_1-t_{0})},
 \quad 0 \;\mathrm{descendant,}
 \\
 {\cal W} \times -p^{-1}(3 f_1^2(x) +    \partial_{x} f_1)(X_\tau)\beta^{-1}e^{\beta(t_1-t_{0})},
 \quad 1\; \mathrm{descendant, type} \; (0)
 \\
 {\cal W} \times -p^{-1}f_1(X_\tau)\beta^{-1}e^{\beta(t_1-t_{0})},
 \quad 1\; \mathrm{descendant, type} \; (1)
 \\
 {\cal W} \times (-3/p)f_1(X_\tau)\beta^{-1}e^{\beta(t_1-t_{0})},
 \quad 2\; \mathrm{descendants, type} \; (0)
 \\
 {\cal W} \times -p^{-1}\beta^{-1}e^{\beta(t_1-t_{0})},
 \quad 2\; \mathrm{descendants, type} \; (0) \; \mathrm{and} \; (1)
 \\
 {\cal W} \times -p^{-1}\beta^{-1}e^{\beta(t_1-t_{0})}, \quad 3\; \mathrm{descendants,}
 \end{array}
 \right.
 \eeaa
 where $0$ descendant means that the particle dies. The particle of type (1) is then simulated on the light cone, meaning that
$Z_{\tau}=\tau$ or   $Z_{\tau}=- \tau$ with probability $1/2$. In the last case, the weight is multiplied by $-1$.
\item For each particle, we apply independently Steps 2 and 3 until the default time -- say $\tau_n$ -- is greater than the maturity $t$. In this case, we multiply $  {\cal W}$ by
 \beaa
 {\cal W} &:=&{\cal W} \times  e^{\beta(t-\tau_{n-1})}
 \eeaa
\item Finally, compute for all particles in $\Kc_t$
 \beaa
 {\cal W} \prod_{k\in\Kc_t} \left( f_2(X^k_t)1_\mathrm{type=0} - f_2(X^k_t-\Delta X^k_t)1_\mathrm{type=1} \right)
 \eeaa
where $\Delta X^k_t=\pm (t-\tau_{n-1})$ with probability $1/2$ and average the result using $M$ Monte-Carlo paths.
\end{enumerate}

\no Below, we have plotted our numerical result for $u(t,X_0)$ again our analytical solution for $d=1$ as a function of $x_0 \in [3,5]$ and $t=1$, see Figure \ref{Fig3}.  We obtain a perfect match.  For completeness, we have also shown the (analytical) solution of the linear wave equation (denoted ``Linearwave'')
 \bea
 (\partial_{tt} - \Delta) u=0, \quad u(0,x)=f_1(x), \quad  \partial_t u(0,x)=f_2(x) \label{nlKG}
 \eea
to show that our numerical solutions capture perfectly the additional nonlinearity $u^3+u \partial_x u$. Note that as $X_t  \in [x_0-t,x_0+t]$, the singularity of $f_2$ at $x=1$ is not relevant.

\begin{figure}[h]
\begin{center}
\includegraphics[width=8cm,height=6cm]{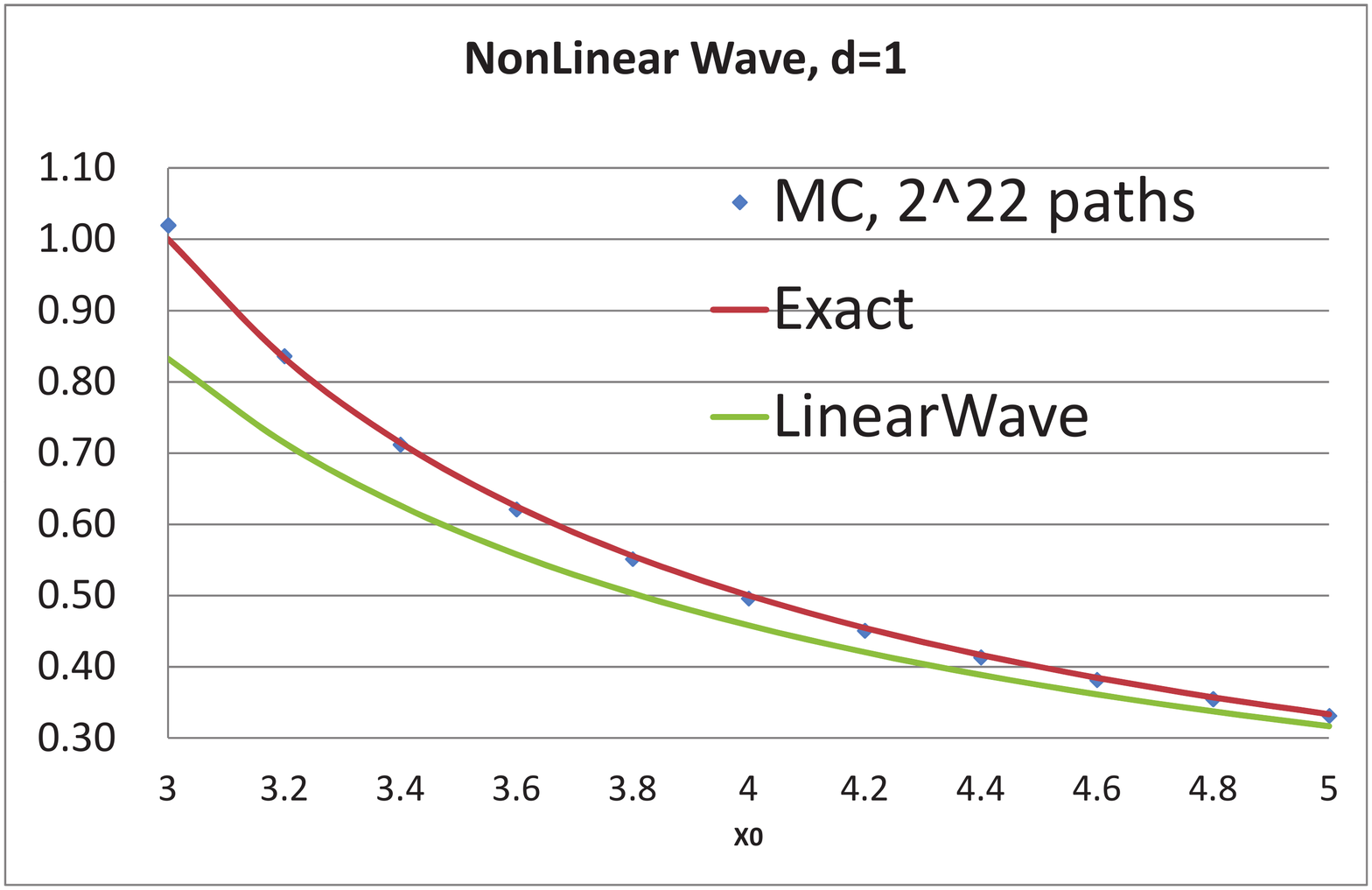}
\end{center}
\caption{Numerical solution of the 1-d semi-linear wave PDE (\ref{scalarYM}) as a function of $x_0$.}
\label{Fig3}
\end{figure}

\subsection{Nonlinear Beam PDE}

We consider the nonlinear beam equation:
 \bea\label{nlbeam}
 \partial_{t}^2 u + \partial_{x}^4 u+u^2 + h(t,x)=0, \quad u(0,x)=\tanh(x), \quad  \partial_t u(0,x)=\cosh(x)^{-2},
 \eea
for which the explicit solution is  $u(t,x)=\tanh( x+ t)$ for a suitable choice of $h$. Here we follow the same discussion as in Section \ref{Wave semi-linear PDE}. Below, we have plotted our numerical result for $u(t,X_0)$ again our exact solution for $d=1$ as a function of $x_0 \in [-0.5,0.5]$ and $t=0.5$, see Figure \ref{Fig3}.  Here we obtain a small error due to the fact that the one-dimensional density $|G(z)| \over \int_\RR |G(z)|dz$  with $G$ given by (\ref{Gbeam}) has been splined on an interval $[-10,10]$ for computational purpose.

\begin{figure}[h]
\begin{center}
\includegraphics[width=8cm,height=6cm]{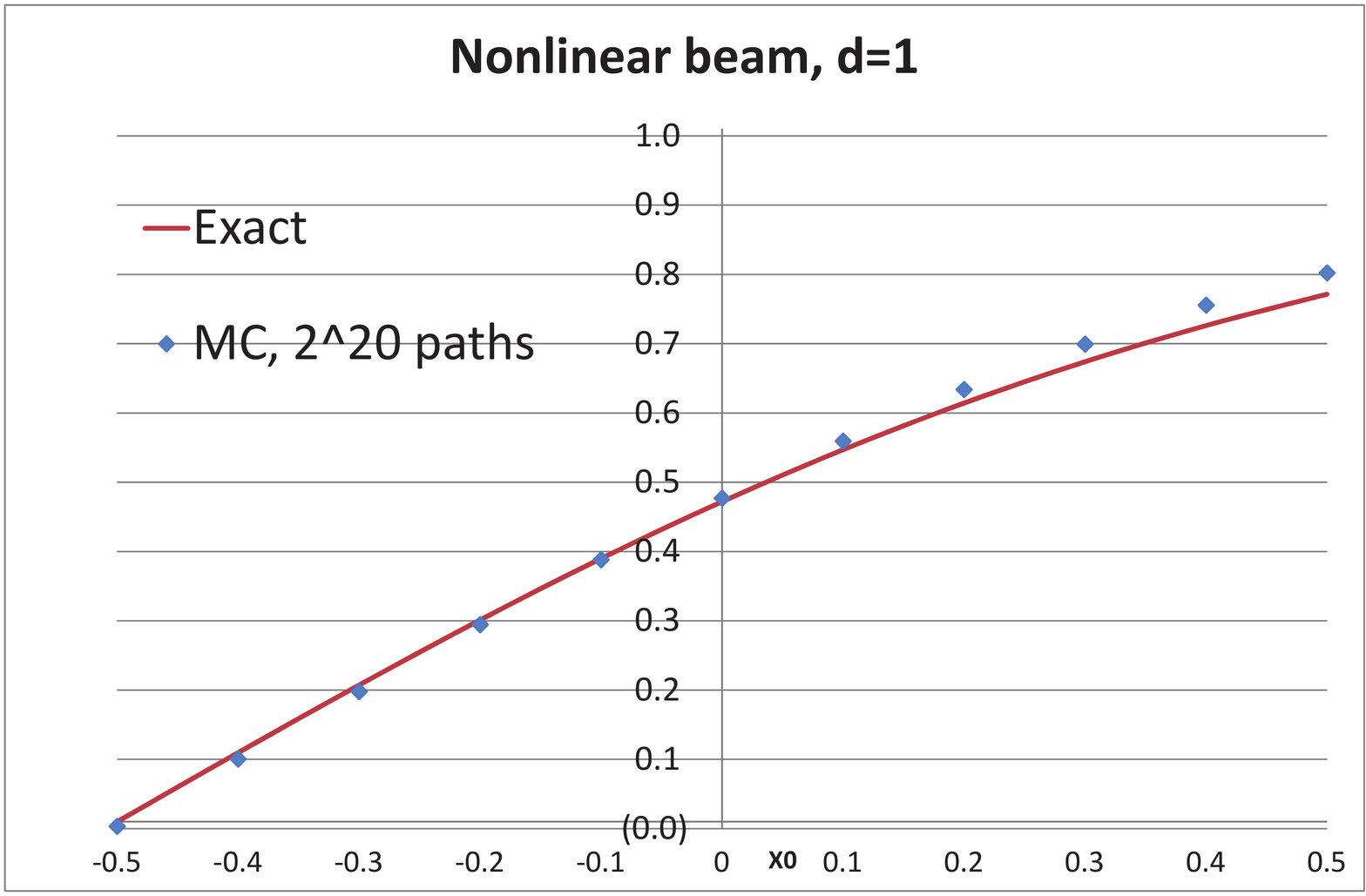}
\end{center}
\caption{Numerical solution of nonlinear beam PDE (\ref{nlbeam}) for $d=1$ as a function of $x_0$.}
\label{Fig3}
\end{figure}

\subsection{Gross-Pitaevskii PDE} \label{Gross--Pitaevskii}
The Gross-Pitaevskii PDE reads
\bea i\partial_t u(t,x)=-{1 \over 2}\Delta u(t,x) +h |u(t,x)|^2 u(t,x), \quad x \in \RR^d \label{GPeq} \eea with $h$ a constant. This equation describes a Bose-Einstein condensate at zero or very low temperature. This has been recently solved using a time-splitting
spectral method \cite{bao}. This deterministic method suffers from the curse of dimensionality and requires suitable mesh size controls. Below, we present our Monte-Carlo algorithm.  We should again emphasize that our algorithm is immune to the dimension and if the standard deviation of our Monte-Carlo estimate converges to zero, we can ensure that we have converged to the true solution.

\no We set $h:=-1$ and  $f_1(x):= {\sqrt{d} \over \cosh(  \sum_{i=1}^d x_i)} $ for which the explicit solution is
 \beaa u(t,x)=e^{i d t \over 2} {\sqrt{d} \over \cosh(  \sum_{i=1}^d x_i)} \eeaa
\no PDE (\ref{GPeq}) can be written as a two-dimensional PDE system with polynomial nonlinearity:
\beaa i\partial_t u(t,x)&=&-{1 \over 2}\Delta u +g u u^* u
\\  -i\partial_t u^*(t,x)&=&-{1 \over 2}\Delta u^* +g u u^* u^* \eeaa  From Example \ref{Schrodinger equation and analytical continuation}, $Z_\tau=e^{i \pi \over 4}\sqrt{\tau} Z \in \mathbb{C}^d$ and $Z^*_\tau=e^{-{i \pi \over 4}}\sqrt{\tau} Z \in \mathbb{C}^d$  with $Z \in N(0,d)$.

\subsubsection*{Meta-algorithm:}
\begin{enumerate}
\item Start at $t_0=0$ at the position $X_0$ with a particle of type $0$ and initialize a {\it complex} weight ${\cal W}:=1$.
\item Simulate an exponential r.v. $\tau$ with (arbitrarily) constant intensity $\beta$ and simulate the particle at the new (complex) position
$X_\tau=X_0+e^{i \pi \over 4}\sqrt{ \tau} Z$ (resp.  $X_\tau=X_0+e^{-{i \pi \over 4}}\sqrt{ \tau} Z$) at $\tau$ if the particle is of type $0$ (resp. of type $1$) with $Z$ a $d$-dimensional standard (real) Gaussian variable.

\no At $t_1:=\tau<t$, the particle dies and we create $2$ descendants of type $0$ and one of type $1$ (resp. $2$ descendants of type $1$ and one of type $0$) if the particle is of type $0$ (resp. type $1$). We then  multiply the weight ${\cal W}$ by
 \beaa
 {\cal W}
 &:=&
 \left\{\begin{array}{l}
 {\cal W} \times (i)\beta^{-1}e^{\beta(t_1-t_{0})},
 \quad \mathrm{type}  \;0
 \\
 {\cal W} \times (-i)\beta^{-1}e^{\beta(t_1-t_{0})} ,
 \quad \mathrm{type}  \;1
 \end{array}
 \right.
 \eeaa
\item For each descendant, we apply independently Steps 2 and 3 until the default time -- say $\tau_n$ -- is greater than the maturity $t$.  In this case, we multiply ${\cal W}$ by
 \beaa
 {\cal W}
 &:=&
 {\cal W} \times  e^{\beta(t-\tau_{n-1})}
 \eeaa
\item Finally, for all particles alive at time $t$ (with locations $(X_t^k)_{k\in\Kc_t}$), compute
 \beaa
 {\cal W} \prod_{k\in\Kc_t \;:\; \mathrm{Type}=0} f_1(X^k_t)\prod_{k\in\Kc_t \;:\; \mathrm{Type}=1} f_1(X^k_t)^*,
 \eeaa
and average the result using $M$ Monte-Carlo paths.
\end{enumerate}

\begin{figure}[h]
\begin{center}
\includegraphics[width=6cm,height=5cm]{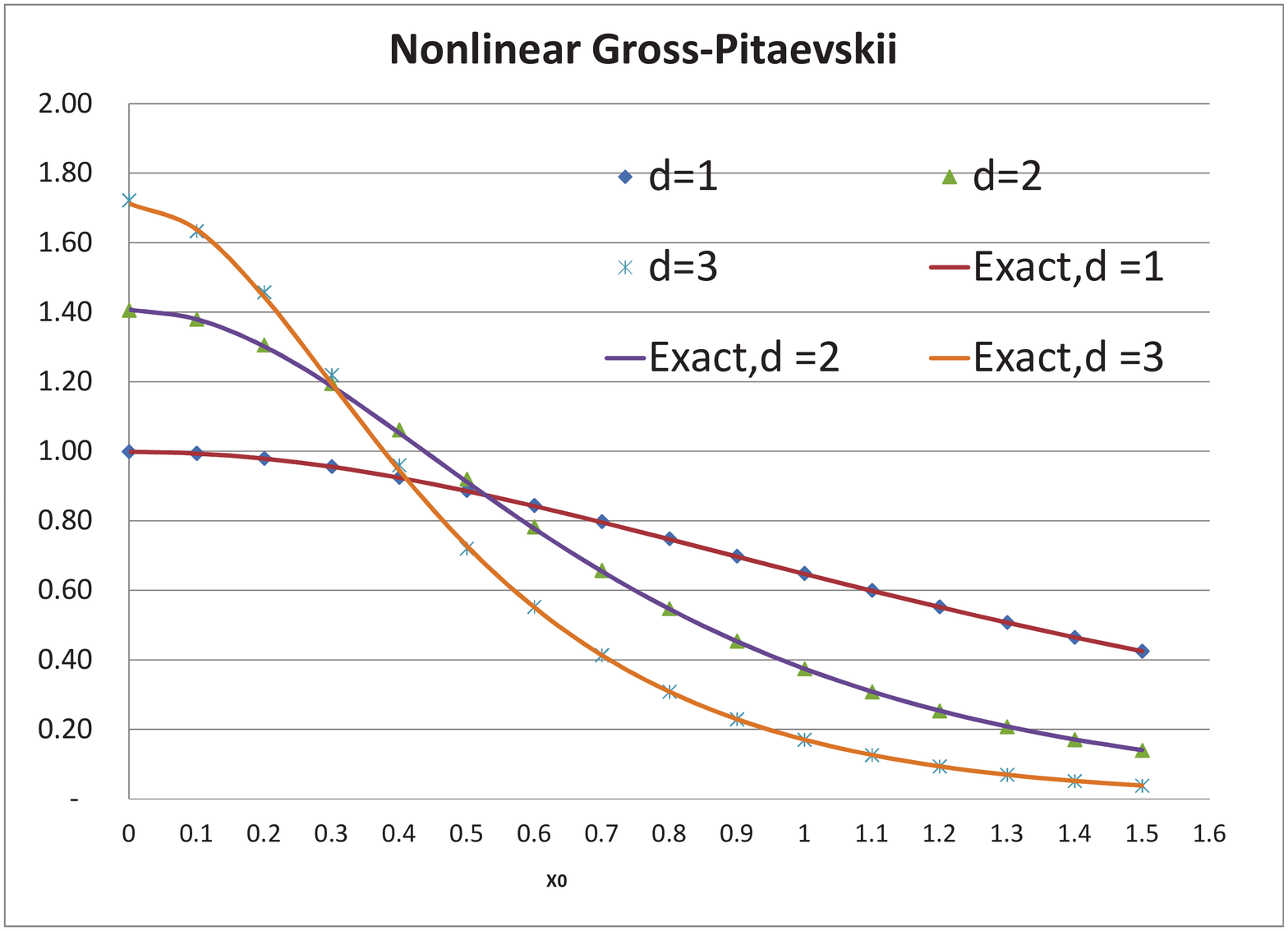}
\includegraphics[width=6cm,height=5cm]{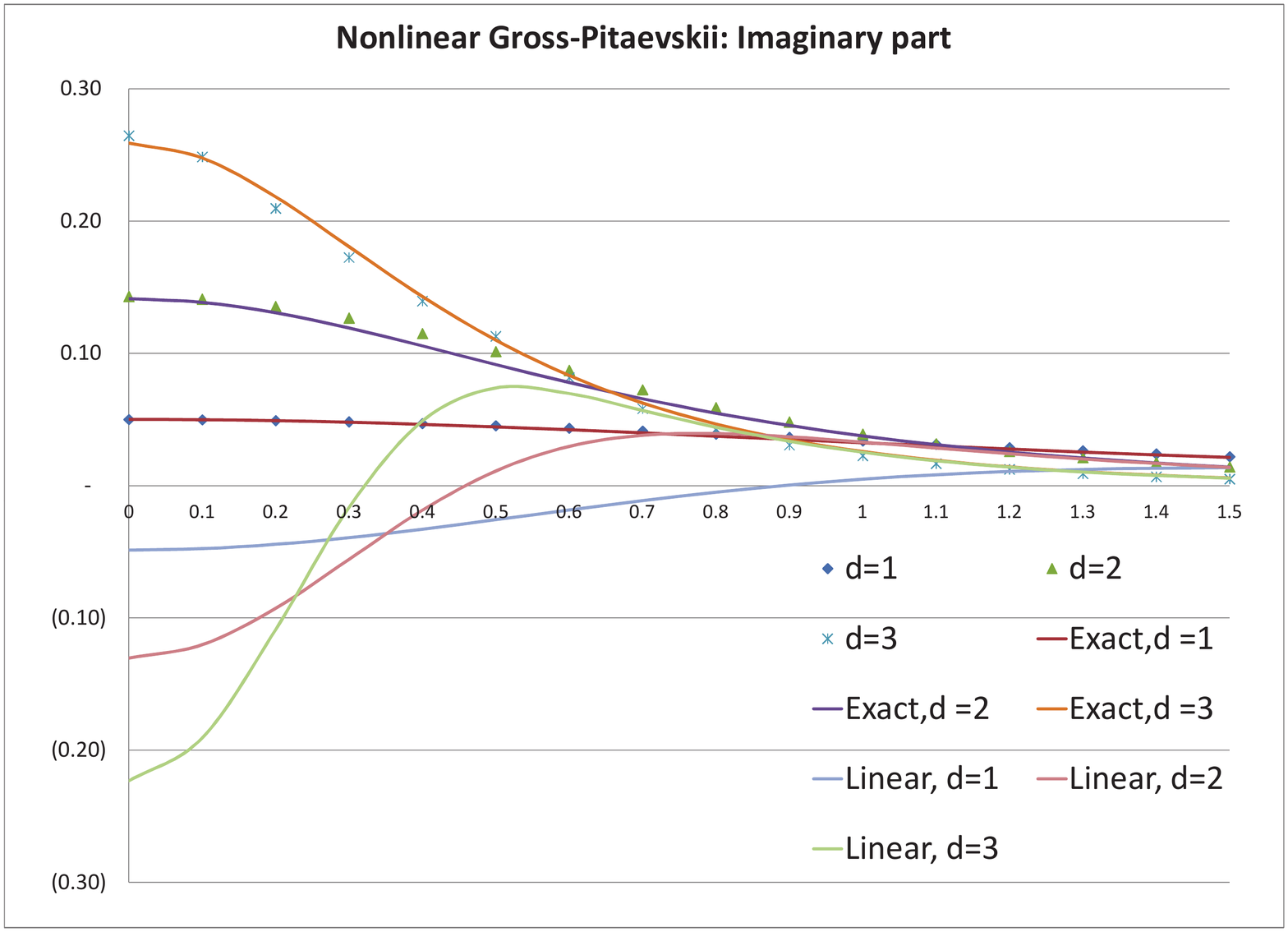}
\end{center}
\caption{Numerical solutions for $Re(u(t,x_0))$ (left) and $Im(u(t,x_0))$ (right) of the Gross-Pitaevskii PDE (\ref{GPeq}) for $d=1,2,3$ as a function of $x_0$.}
\label{FigGP}
\end{figure}

\no We have chosen $\beta=1$ in our numerical experiments and $M=2^{22}$ for which the standard deviation of our estimator is less than $0.01$. Below, we have plotted our numerical result for $Re(u(t,x_0))$  and $Im(u(t,x_0))$ again our analytical solution for $d=1,2,3$ as a function of $x_0 \in [0,1.5]$ (all the coordinates in $\RR^d$ are equal to $x_0$) and $t=0.1$, see Figure \ref{FigGP}.  We obtain a perfect match.  For completeness, we have also shown the (numerical) solution of the linear Schrodinger equation (i.e., $g:=0$),  denoted ``Linear''.

\end{document}